\documentclass[9pt,final]{IEEEtran}

\usepackage[cmky]{xcolor}

\usepackage{caption}
\usepackage{graphicx}

  \usepackage{rotating}
  \usepackage{floatpag}
  \rotfloatpagestyle{empty}

  \usepackage{amsmath,amssymb}  
\usepackage{wrapfig}

\usepackage{mathrsfs,multicol}
\usepackage[mathscr]{euscript}

\newtheorem{theorem}{Theorem}
\newtheorem{corollary}{Corollary}[theorem]
\newtheorem{lemma}
{Lemma}
 
\newtheorem{definition}
{Definition}
\newtheorem{problem}
{Problem}
{Example}

\usepackage{tikz}
\usepackage{pgfplots}
\usepackage{tikz-3dplot}

\usetikzlibrary{decorations.pathmorphing} 
\usetikzlibrary{matrix} 
\usetikzlibrary{arrows.meta} 
\usetikzlibrary{shapes} 
\usetikzlibrary{calc} 
\usetikzlibrary{positioning} 
\usetikzlibrary{patterns}
\usetikzlibrary{intersections}
\usepgfplotslibrary{fillbetween}
\usepgfplotslibrary{polar}
\usepgfplotslibrary{colorbrewer}

\usepgfplotslibrary{groupplots}
\usepgfplotslibrary{fillbetween}

\usetikzlibrary{spy}

\tikzstyle{block} = [draw,rectangle,thick,minimum height=2em,minimum width=2em]
\tikzstyle{sum} = [draw,circle,inner sep=0mm,minimum size=4mm]
\tikzstyle{connector} = [->,thick]
\tikzstyle{line} = [thick]
\tikzstyle{branch} = [circle,inner sep=0pt,minimum size=1mm,fill=black,draw=black]
\tikzstyle{guide} = []
\tikzstyle{legendBlock} = [rectangle,minimum height=2em,minimum width=2em]

\makeatletter    
\pgfplotsset{    
    compat=1.13,
    /tikz/max node/.style={
        anchor=south,
    },
    /tikz/min node/.style={
        anchor=south,
        name=minimum
    },
    mark min/.style={
        point meta rel=per plot,
        scatter/@pre marker code/.code={%
            \ifx\pgfplotspointmeta\pgfplots@metamin
                \def\markopts{}%
                \coordinate (minimum);
                \node [min node] {
                };
            \else
                \def\markopts{mark=none}
            \fi
            \expandafter\scope\expandafter[\markopts,every node near
coord/.style=green]
        },%
        scatter/@post marker code/.code={%
            \endscope
        },
        scatter,
    },
    mark max/.style={
        point meta rel=per plot,
        visualization depends on={x \as \xvalue},
        scatter/@pre marker code/.code={%
        \ifx\pgfplotspointmeta\pgfplots@metamax
            \def\markopts{}%
            \coordinate (maximum);
            \node [max node] {
                \pgfmathprintnumber[fixed]{\xvalue},%
                \pgfmathprintnumber[fixed]{\pgfplotspointmeta}
            };
        \else
            \def\markopts{mark=none}
        \fi
            \expandafter\scope\expandafter[\markopts]
        },%
        scatter/@post marker code/.code={%
            \endscope
        },
        scatter
    }
}
\makeatother    

\definecolor{yellow1}{RGB}{255,255,204}
\definecolor{blue2}{RGB}{161,218,180}
\definecolor{blue3}{RGB}{65,182,196}
\definecolor{blue4}{RGB}{44,127,184}
\definecolor{blue5}{RGB}{37,52,148}
\definecolor{brown1}{RGB}{166,97,26}
\definecolor{brown2}{RGB}{223,194,125}
\definecolor{red1}{RGB}{215,25,28}
\definecolor{red2}{RGB}{253,174,97}

\definecolor{blackOne}{cmyk}{0,0,0,85}
\definecolor{grayThree}{cmyk}{0,0,0,41}

\definecolor{darkViolet}{cmyk}{65,70,0,0}
\definecolor{darkSea}{cmyk}{85,30,0,0}
\definecolor{paleSea}{cmyk}{33,3,0,0}
\definecolor{paleViolet}{cmyk}{65,70,0,0}
\definecolor{paleGreen}{cmyk}{24,0,39,0}
\definecolor{paleOrange}{cmyk}{5,35,70,0}

\definecolor{brick-red-con}{RGB}{222,66,91}
\definecolor{orange-con}{RGB}{225,164,124}
\definecolor{cyan-con}{RGB}{72,143,49}
\definecolor{light-green-con}{RGB}{161,185,149}

\definecolor{redSeries1}{RGB}{254,240,217}
\definecolor{redSeries2}{RGB}{253,204,138}
\definecolor{redSeries3}{RGB}{252,141,89}
\definecolor{redSeries4}{RGB}{215,48,31}

\definecolor{tolYellowMain}{HTML}{DDAA33}
\definecolor{tolRedMain}{HTML}{BB5566}
\definecolor{tolBlueMain}{HTML}{004488}
\definecolor{tolPaleYellow}{HTML}{EEEEBB}
\definecolor{tolPaleRed}{HTML}{FFCCCC}
\definecolor{tolDarkRed}{HTML}{663333}

\usepackage{algpseudocode,algorithm,algorithmicx}

\algrenewcommand\algorithmicrequire{\textbf{Precondition:}}
\algrenewcommand\algorithmicensure{\textbf{Postcondition:}}
\algnewcommand{\Inputs}[1]{%
      \State \textbf{Inputs:}
        \Statex \hspace*{\algorithmicindent}\parbox[t]{.8\linewidth}{\raggedright #1}
    }

\usepackage{subfig}

\pgfplotstableread{
    xcoord       ycoord
             0 1.35  
             0.1 1.15  
             0.25 0.35  
             0.5 0.1  
             0.55 0  
             0.5 -0.1  
             0.25 -0.35  
             0.1 -1.15  
             0 -1.35  
             -0.1 -1.15  
             -0.25 -0.35  
             -0.5 -0.1  
             -0.55 0  
             -0.5 0.1  
             -0.25 0.35  
             -0.1 1.15  
             0 1.35  
             0 1.35  
             0.1 1.15  
             0.25 0.35  
             0.5 0.1  
             0.55 0  
             0.5 -0.1  
             0.25 -0.35  
             0.1 -1.15  
             0 -1.35  
             -0.1 -1.15  
             -0.25 -0.35  
             -0.5 -0.1  
             -0.55 0  
             -0.5 0.1  
             -0.25 0.35  
             -0.1 1.15  
             0 1.35  
}\data

\usepackage{mathtools}
\usepackage{amsfonts}
\usepackage[T1]{fontenc}
\usepackage{ae,aecompl}

\begin{document}
\title{Estimating a scalar log-concave random variable, using a silence set based probabilistic sampling}
\author{Maben~Rabi, 
        Junfeng~Wu, %
        Vyoma~Singh, %
        Karl~Henrik~Johansson%
\thanks{M.~Rabi is with
{\O}stfold University College, Norway. e-mail: firstname.lastname@hiof.no.
J.~Wu is with Chinese University of Hong Kong, Shenzen, China.
V.~Singh is with Koneru Lakshmaiah Education Foundation, Guntur, India.
K.~H.~Johansson is with the Royal Institute of technology~(KTH), Sweden.}
}

\maketitle

\begin{abstract}%
    We study 
    the probabilistic sampling of a random variable, in which the variable is sampled only if it falls outside a given set, which is called the silence set.
    This helps us to understand optimal event-based sampling for the special case of IID random processes, and also to understand the design of a sub-optimal scheme for other cases.
    We consider the design of this probabilistic sampling for a scalar, log-concave random variable, to minimize  either the mean square estimation error, or the mean absolute estimation error. We show that the optimal silence interval: (i)~is essentially unique, and (ii)~is the limit of an iterative procedure of centering. Further we show through numerical experiments that super-level intervals seem to be remarkably near-optimal for mean square estimation. Finally we use the Gauss inequality for scalar unimodal densities, to show that probabilistic sampling gives a mean square distortion that is less than a third of the distortion incurred by periodic sampling, if the average sampling rate is between 0.3~and 0.9~samples per tick.
\end{abstract}


\section{Introduction}
Consider a sensor and  a tracking station, that are connected by an ideal, analog communication link from the sensor to the tracking station. The sensor takes perfect observations of a discrete-time random 
process. 
At times of its choice, the sensor sends its current samples.
The sensor is allowed to choose these times on the fly, based upon the causal record of its transmission decisions and transmitted messages up till the previous time instant.
At times when the sensor does not send samples, it sends a special {\textsf{SILENCE}} symbol.
This scheme is called {\textit{Event-based sampling~\cite{miskowicz2015eventBasedControlBook},}} 
and is worthwhile if:
\begin{itemize}
\item{sending and reliably receiving the {\textsf{SILENCE}} symbol costs a negligible expenditure of energy, or bandwidth, and}
\item{sending and reliably receiving the observed measurement costs a significant expenditure of energy, or bandwidth.} 
\end{itemize}
To immplement this form of sampling, we need to causally prescribe a {\textit{silence set}} for each instant of time. Samples are generated or not, dependent on whether or not the observed signal falls outside the silence sets.


\subsection{Previous results}
Broadly speaking, two approaches can be seen in the literature -  a deterministic one~\cite{lemmon2009sienaSurvey,heemelsJohanssonTabuada2012}, and a stochastic one~\cite{astrom-bernhardsson-2003,shiShiChen2015eventBasedEstimationBook}.
In the deterministic approach, silence set design is posed as a problem of ensuring that a candidate Lyapunov function decreases over time, or a suitable objective function of the state and control signals is minimized. In the stochastic approach, silence set design is posed as a networked sequential design problem~\cite{rabiRameshJohansson2016siam,jhelumChakravorthyMahajan2020remoteEstimation}, with two or more decision agents. Within the stochastic approach, a variety of communication limitations have been captured: limited packet rate~\cite{jhelumChakravorty2017fundamentalLimits}, packet losses~\cite{jhelumChakravorthyMahajan2020remoteEstimation} and delays~\cite{khojasteh2020valueOfTimingInformation}.   Both approaches are computationally demanding, as they  require the solution of Linear or Bilinear matrix inequalities, or a Dynamic programming problem.

If: (i)~the state signal is scalar, (ii)~the dynamics is linear,
and (iii)~the process and sensor noise densities are Gaussian, then the estimation error variance is minimized by symmetric silence intervals around the Kalman predictor~\cite{hajekPagingJournal,lipsaMartins2011,nayyarAndThreeProfessorsTransactions2013,jhelumChakravorthyMahajan2020remoteEstimation}. With this simplified structure for the optimal silence sets, the calculation of their sizes requires the numerical solution of an one-agent Dynamic programming problem.  
Andr{\'e}n~et~al.~\cite{andren2017cdcPaperWithNonConvexSilenceSets} have shown that optimal silence sets can be non-convex.

 Molin and Hirche~\cite{molin2012unimodalBimodal,molinHirche2017iterativeSilenceSets} give two-person iterative algorithm that they show to globally converge to a sequence of symmetric intervals around the Kalman predictor, if the initial state and noise are unimodal and symmetric, as assumed in the works mentioned in the previous paragraph. 

A log-concave probability density~\cite{an1996logConcave} is one whose logarithm is a concave function over its support.
Henningsson and {\AA}str{\"o}m~\cite{henningssonAstrom2006mtnsLogConcave} showed that if the noise densities are  log-concave, then optimal observer design for linear systems is simpler than without log-concavity, because the density of the estimation error is log-concave even after conditioning on the signal being inside the  silence set.  For s special class of log-concave densities, called strongly log-concave, in which tails decay at least as fast as for Gaussian densities, they give an upper bound on the variance of estimation error.
\subsection{Ultra-myopic choice of silence sets}
To sidestep the computational burden of calculating optimal silence sets, we use the suboptimal strategy of designing each individual silence set without regard for the effects of this choice on the costs and constraints at later times instants.  Given a causal specification of a sampling rate budget, our strategy is to compute a suboptimal silence set for each time such that 
it minimizes only the estimation distortion at that instant, with no concern for the  distortion at later time instants. 
\subsection{Our results}
We pursue the question of how best to 
synthesize this suboptimal strategy, when the observed signal is a scalar, and is driven by noises that have log-concave densities.

In Section~\ref{section:singleStageVector} we give the centering algorithm  of ~\cite{molin2012unimodalBimodal,molinHirche2017iterativeSilenceSets} for improving a given silence set.  We are able to prove  in Section~\ref{section:scalarlog-concave} that for any scalar log-concave density,  repeated runs of the centering  algorithm converge to a unique optimal silence interval. 
This  is similar to 
the result~\cite{trushkin1982,trushkin1984,kieffer1983uniquenessOfOptimalQuantizer} that every scalar log-concave density has  a unique optimal quantizer.

When a scalar log-concave density is symmetric, then clearly it is optimal to use   a symmetric silence interval around the mean cum mode. When such a density is not symmetric, we posit that super-level sets give  excellent performance. We give some empirical evidence in Section~\ref{section:superLevelInterval} to support this claim.

Finally, in Section~\ref{section:rateDistortionBound}
we use the Gauss inequality for unimodal densities to bound the rate-distortion trade-off incurred by using silence intervals that are symmetric around the mode.

\section{Centering reduces distortion\label{section:singleStageVector}}

We start with a general formulation of our problem, and specialize to the scalar case in Section~\ref{section:scalarlog-concave}.
\begin{problem}
For 
the random vector~$X \in  {\mathbb{R}}^n ,$ the probability density 
exists, and is 
given. The probabilistic sampling problem 
is to pick a measurable silence set~${\mathscr{A}},$ 
by an  optimization in which: 
(a)~the following chance constraint is met:
\begin{align}
 {\mathbb{P}} \left[X \in {\mathscr{A}}  \right] & \ge \eta, \ {\text{for some prescribed}} \ \eta \in 
     \left[ 0, 1\right], \label{eqn:chanceConstraintVector}
\end{align}
and (b)~the set~${\mathscr{A}}$ minimizes the average distortion:
\begin{gather*}
    {\mathbb{E}} \left[ \delta {\left( X - 
          {\widehat{X}}_{\mathscr{A}}
                     \right)} 
\left| X \in {\mathscr{A}} \right. 
                     \right],
\end{gather*}
where~$ {\widehat{X}}_A $ denotes the estimate under silence, and is defined as:
\begin{align} 
     {\widehat{X}}_{\mathscr{A}}
     & \triangleq 
    \arg\min_{Y\in{\mathbb{R}}^n}{
        {\mathbb{E}} \left[\delta\left( X - Y \right)
                            \left \lvert X\in {\mathscr{A}} \right.
                     \right] 
        } ,
\label{eqn:bestEstimateDefinition}
\end{align}
and where the distortion function~$ \delta : 
  {\mathbb{R}}^n \to {\mathbb{R}}
$
is 
non-negative, and
monotonically increasing as we move radially away from the origin along an arbitrary vector. 
\label{pbm:optimalSilenceSet}
\end{problem}
In other words, for any vector~$ x \in {\mathbb{R}}^n , $ the value~$\delta\left( \alpha x \right)$ is an increasing function of the positive real number~$\alpha.$ Hence, given two different positive levels, the sub-level set corresponding to the smaller level is guaranteed to be entirely within the interior of the sub-level set corresponding to the higher level.
Examples of valid distortion functions in the variable~$x$ are:
\begin{gather*}
   {\left\lvert x \right\rvert},  \, {\left\lvert x \right\rvert}^2,  \ \text{or} \
    {x }^T M {x} \; \text{with the matrix} \; M > 0.  
\end{gather*}

\subsection{Centering is necessary for optimality}

\begin{definition}[centering of a set]
     Suppose we are given  
     the random variable~$X, $   
     the silence set~$ {\mathscr{S}} ,$ 
     the probability value~$\eta,$ the distortion function~$ \delta ( \cdot ) ,$ 
     and the best estimate~$ {\widehat{X}}_ {\mathscr{S}}  $ 
     as definied by Equation~\ref{eqn:bestEstimateDefinition}.
     Then the centering of the set~${\mathscr{S}}$ is denoted
by~$ {\mathscr{S}}_{\text{cen.}} , $  and is defined as the smallest sub-level set of the function~%
$\delta\left(  \boldsymbol{\cdot} -  {\widehat{X}}_ {\mathscr{S}} 
      \right)
,
$
 such that the set has a probability mass no less than~$ \eta . $
    In other words, if we let the nonegative number
\begin{align}
    \alpha_{\eta} &  \triangleq \inf{ \left\{ \alpha : 
                        {\mathbb{P}}\left[ 
                        X  \left\lvert  
                             \delta \left( X - 
                    {\widehat{X}}_{\mathscr{S}}
                             \right) \le \alpha 
                          \right.
                      \right]  \ge \eta \right\}
                      },  \ \text{then}
                      \nonumber \\
    {\mathscr{S}}_{\text{cen.}}  & = 
                      \left\{
                              x \in {\mathbb{R}}^n 
                        :
                             \delta \left( x - 
                    {\widehat{X}}_{\mathscr{S}}
                             \right) \le \alpha_{\eta} 
                          \right\}
                    .
\end{align}
\end{definition}

For any two sets ${\mathscr{A}}, {\mathscr{B}}, $ we denote by ${\mathscr{A}}/{\mathscr{B}}$ the set ${\mathscr{A}} \cap {\mathscr{B}}^C,$ which is the set of all those elements of~${\mathscr{A}}$ not found in~${\mathscr{B}}.$

A {\textit{centered set}} is any set that either equals its centering, or differs from it by a set of probability zero.  
This means that:
 \begin{align*}
      {\mathbb{P}} \left[ X \in {\mathscr{S}} /  {\mathscr{S}}_{\text{cen.}}
                     \right] 
                     &
                     = 0 ,
 \ \ 
 \text{if and only if~$ {\mathscr{S}} $ is a centered set.}
 \end{align*}
   For example consider sets on the real line, and let~$\delta\left(\cdot\right)$ be the usual squared error distortion. If an interval has its conditional mean equal to its midpoint, then this interval is centered.
\begin{algorithm}
\caption{Centering algorithm%
\label{alg:centering}%
}
\begin{algorithmic}[1]
     \Statex
     \Function{Centering}{$ \ p_X{\left(\cdot\right)}, \ {\mathscr{S}}, \ \delta\left( \cdot \right) , \ \eta $}
     \State{$
     {\widehat{X}}_{\mathscr{S}}
      =
    \arg\min_{Y\in{\mathbb{R}}^n}{
        {\mathbb{E}} \left[\delta\left( X - Y \right)
                            \left \lvert X\in {\mathscr{S}} \right.
                     \right] 
        } $    
        }
\State{$ \alpha_{\eta}  
           = \inf{\left\{ 
                     \alpha \in {\mathbb{R}} 
              :
                      {\mathbb{P}}\left[ 
                           X 
                          \left\lvert  
                             \delta \left( x - 
                 {\widehat{X}} \right) \le \alpha  
                                 \right]
                     \ge \eta
                  \right. \;  \right\}}
           $}
\State{${\mathscr{S}}   =
                          \left\{
                              x \in {\mathbb{R}}^n 
                    :  \delta \left( x -   {\widehat{X}} \right) \le \alpha_{\eta}  
                          \right\}
                      $}
    \Comment{\textcolor{blue4}{$ {\mathscr{S}} \longleftarrow {\mathscr{S}}_{\text{cen.}}  $}}
\State\Return{${\mathscr{S}} $}
    \EndFunction
\end{algorithmic}
\end{algorithm}

\begin{lemma}
Let the random vector~$X$ have a regular probability density function. Given the probability value~$ \eta , $  and the distortion function~$ \delta\left(\cdot\right) ,$ consider  Problem~\ref{pbm:optimalSilenceSet}. If an optimum silence set exists, then a centered one 
exists that is optimal.
\label{example:circularDisk}
\end{lemma}
\begin{proof}
We shall show that if a candidate optimal set is such that it differs from its centering by a set of non-zero probability, then centering lowers the distortion.

Let the random vector~$X\in{\mathbb{R}}^n . $ 
Suppose that the positive number $ \eta $ is such that~$\eta\in(0,1)$.  Then consider any candidate optimal set~${\mathscr{A}}^*$ that is not centered, and satisfies the 
constraint:
\begin{align*}
    {\mathbb{P}} \left[  
    {\mathscr{A}}^* \right] & \ge \eta . 
\end{align*}
Let~$  {\widehat{X}}_{ {\mathscr{A}}^* }  $ denote
the best estimate under the silence set~$  {\mathscr{A}}^*  . $

Without loss of generality, we can assume that the candidate optimal set has a probability mass of exactly~$\eta .$ If this were not the case, then we can shrink the given set, to derive a new silence set that has a probability mass of exactly~$ \eta ,$ but with a lower average distortion. In specific, given the candidate set, and a shrinkage factor~$\sigma,$ with~$ 0 < \sigma < 1 ,$ the shrunken set is:
\begin{gather*}
    \left\{  
        x \in  {\mathbb{R}}^n : 
        {\frac{1}{\sigma}} \cdot 
           {\left( x - {\widehat{X}}_{ {\mathscr{A}}^* } \right) }
           \in  {\mathscr{A}}^*
    \right\} .
\end{gather*}
Because the density is regular, the shrinkage factor~$ \sigma $ can be chosen to achieve a probability mass of exactly~$ \eta $ for the shrunken set. Hence, we shall assume that the given candidate optimal silence set~$  {\mathscr{A}}^*  $ satisfies the exact chance equality:
\begin{align}
    {\mathbb{P}} \left[  
    {\mathscr{A}}^* \right] & = \eta . %
\end{align}

Then consider 
the centering of the given set~$  {\mathscr{A}}^* $:
\begin{align*}
    C_{\eta} & \triangleq \left\{ x \in {\mathbb{R}}^n : 
      \delta{\left( x - {\widehat{X}}_{ {\mathscr{A}}^* } \right)} \le r_{\eta} \right\},
\end{align*}
    where the `radius'~$r_{\eta}$ is chosen such that
$
    {\mathbb{P}} \left[  
                         C_{\eta} \right] = \eta.
$

A suitable value for~$ _{\eta} $ can be chosen, because the density is regular.  
Since
~$\eta$ is strictly less than~$1,$ the radius of~$C_{\eta}$ must be finite. By assumption, the set~$ {\mathscr{A}}^* $ is not centered. Hence:
\begin{align*}
   {\mathbb{P}} \left[  
      {\mathscr{A}}^{\star} / C_{\eta} \right] 
      & > 0 . 
\end{align*}
\pgfmathsetseed{1}
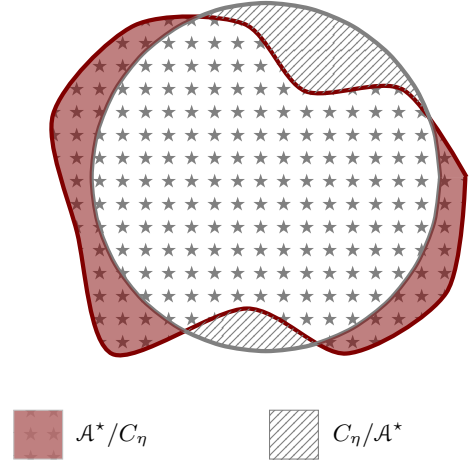
\begin{figure}
\begin{center}
\begin{tikzpicture}[scale=1]
\begin{scope}
\begin{polaraxis}[hide axis]
\addplot [ultra thick, color=red!50!black, domain=0:360, samples=12, smooth, name path = setB, pattern = fivepointed stars, pattern color = gray,]
  {2-rand};
\addplot [ultra thick, color=gray, domain=0:360, samples=25, smooth, name path = circle]
  {2.0};
\tikzfillbetween[ 
    of = circle and setB,
    split, 
    every even segment/.style={fill=red!50!black, opacity=0.5},
    every odd segment/.style={fill,pattern=north east lines,pattern color=gray}, 
  ]{gray};
\end{polaraxis}
\end{scope}
\begin{scope}{xshift=0,yshift=-2cm}
\node[legendBlock,draw=gray,fill=red!50!black, opacity=0.5,
       preaction={  pattern = fivepointed stars, pattern color = gray},
     ] at (0.4,0) {}; 
\node[anchor = west ] at (0.8,0) {${\mathscr{A}}^{\star}/C_{\eta}$}; 
\node[legendBlock, draw = gray, pattern = north east lines, pattern color = gray ] at (3.8,0) {}; 
\node[anchor = west ] at (4.2,0) {$C_{\eta}/{\mathscr{A}}^{\star}$}; 
\end{scope}
\end{tikzpicture}
\end{center}
\caption{The sets used in the proof of Lemma~\ref{example:circularDisk}.}
\label{fig:spilloverAndSpillunderSets}
\end{figure}
Since  both 
$
 {\mathscr{A}}^{\star} 
  ,$
    and
$
    C_{\eta}  
$
have the same probability mass, we get:
\begin{align*}
   {\mathbb{P}} \left[ 
   {\mathscr{A}}^{\star} / C_{\eta} \right] 
& =
   {\mathbb{P}} \left[ 
   C_{\eta} / {\mathscr{A}}^{\star} \right].
\end{align*}
Note that on every point~$x$ of the set~$ {\mathscr{A}}^{\star} / C_{\eta},$
 the distortion function~$ 
  \delta{\left( x -  {\widehat{X}}_{{\mathscr{A}}^{\star}}  \right)}
$ 
takes values that are bigger than 
values 
of the function 
at points of the set~$ C_{\eta} / {\mathscr{A}}^{\star}.$
Therefore, 
\begin{align*}
 {\mathbb{E}} \left[ \left.
\delta{ \left( X -  {\widehat{X}}_{{\mathscr{A}}^{\star}}  \right) }
\right|  X \in {\mathscr{A}}^{\star} / C_{\eta} \right]    
& > 
 {\mathbb{E}} \left[ \left.
\delta{ \left( X -  {\widehat{X}}_{{\mathscr{A}}^{\star}}  \right) }
\right|  X \in  C_{\eta} / {\mathscr{A}}^{\star} \right]
.
\end{align*}
Using the above inequalities, we can write as below (at each step, the text in  {\textcolor{red}{red}} signifies changes from the previous step):
\begin{multline*}
\shoveleft{ \mspace{-28.0mu}
 {\mathbb{E}} \left[ \left.
\delta{ \left( X -  {\widehat{X}}_{{\mathscr{A}}^{\star}} \right) }
\right|  
  {\mathscr{A}}^{\star}  \right] 
 \, \cdot \, {\mathbb{P}} \left[ 
        {\mathscr{A}}^{\star}   \right]
\phantom{\; = \;
\int_{ {\mathscr{A}}^{\star} }{  
  \delta{ \left( x -  {\widehat{X}}_{{\mathscr{A}}^{\star}} \right) }
\, \quad \quad \quad \quad \quad 
}}
}
\\
\begin{split}
& = \
\int_{ {\mathscr{A}}^{\star} }{  
  \delta{ \left( x -  {\widehat{X}}_{{\mathscr{A}}^{\star}} \right) }
\,  p_X\left(x\right) dx
}
\\
& = \
 \int_{ 
     {\textcolor{red}{%
 {\mathscr{A}}^{\star} / C_{\eta} 
   }}
 }{  
  \delta{ \left( x -  {\widehat{X}}_{{\mathscr{A}}^{\star}} \right) }
\,  p_X\left(x\right) dx
}
\\
& \quad \ \ + 
 \int_{ 
      {\textcolor{red}{%
   {\mathscr{A}}^{\star} \cap \, C_{\eta} 
     }}
   }{  
  \delta{ \left( x -  {\widehat{X}}_{{\mathscr{A}}^{\star}} \right) }
\,  p_X\left(x\right) dx
} ,
\\
& = \
 {\mathbb{E}} \left[ \left.
\delta{ \left( X -  {\widehat{X}}_{{\mathscr{A}}^{\star}}  \right) }
\right|  
{\mathscr{A}}^{\star} / C_{\eta} \right] 
 \cdot {\mathbb{P}} \left[ 
    {\mathscr{A}}^{\star} /  C_{\eta}  \right]
 \\
& \quad  \ \ +
 {\mathbb{E}} \left[ \left.
\delta{ \left( X -  {\widehat{X}}_{{\mathscr{A}}^{\star}}  \right) }
\right| 
   {\mathscr{A}}^{\star} \cap C_{\eta} \right] 
 \cdot {\mathbb{P}} \left[ 
     {\mathscr{A}}^{\star} \cap \,  C_{\eta} \right]
\\
& > \
 {\mathbb{E}} \left[ \left.
\delta{ \left( X -  {\widehat{X}}_{{\mathscr{A}}^{\star}}  \right) }
\right| 
  {\textcolor{red}{%
C_{\eta} / {\mathscr{A}}^{\star} 
 }}
\right] 
 \cdot {\mathbb{P}} \left[ 
       {\mathscr{A}}^{\star} /   C_{\eta}   \right]
 \\
& \quad  \ \ +
 {\mathbb{E}} \left[ \left.
\delta{ \left( X -  {\widehat{X}}_{{\mathscr{A}}^{\star}}  \right) }
\right|  
 {\mathscr{A}}^{\star} \cap C_{\eta} \right] 
 \cdot {\mathbb{P}} \left[ 
      {\mathscr{A}}^{\star} \cap  C_{\eta}    \right] ,
\\
& = \
 {\mathbb{E}} \left[ \left.
\delta{ \left( X -  {\widehat{X}}_{{\mathscr{A}}^{\star}}  \right) }
\right|  
   C_{\eta} / {\mathscr{A}}^{\star}  \right] 
 \cdot {\mathbb{P}} \left[ 
     {\textcolor{red}{%
     C_{\eta} /  {\mathscr{A}}^{\star}  
     }}
     \right]
 \\
& \quad  \ \ +
 {\mathbb{E}} \left[ \left.
\delta{ \left( X -  {\widehat{X}}_{{\mathscr{A}}^{\star}}  \right) }
\right|  
{\mathscr{A}}^{\star} \cap C_{\eta} \right] 
 \cdot {\mathbb{P}} \left[ 
      {\mathscr{A}}^{\star} \cap   C_{\eta}   \right] ,
\\
& = \
 {\mathbb{E}} \left[ \left.
 \delta{ \left( X -  {\widehat{X}}_{{\mathscr{A}}^{\star}}  \right) }
\right| 
  C_{\eta} \right] 
 \cdot {\mathbb{P}} \left[ C_{\eta}   \right]
      ,
\\
& \ge \
 {\mathbb{E}} \left[ \left.
 \delta{ \left( X - 
     {\textcolor{red}{%
 {\widehat{X}}_{  C_{\eta} }
     }}
 \right) }
\right|  
   C_{\eta} \right] 
 \cdot {\mathbb{P}} \left[ C_{\eta}   \right] .
\end{split}
\end{multline*}
Since both the sets~$ {{\mathscr{A}}^{\star}}  ,  C_{\eta}   $
have the same probability mass, 
\begin{align*}
{\mathbb{E}} \left[ \left.
 \delta{ \left( X - 
 {\widehat{X}}_{  {{\mathscr{A}}^{\star}}  }
 \right) }
\right|  
   {{\mathscr{A}}^{\star}}  \right] 
& > \
 {\mathbb{E}} \left[ \left.
 \delta{ \left( X - 
 {\widehat{X}}_{  C_{\eta} }
 \right) }
\right|  
   C_{\eta} \right] 
.
\end{align*}
Hence centering cannot increase the average distortion.
\end{proof}

\section{Uniqueness of the optimal silence interval\label{section:scalarlog-concave}}
Each run of the centering algorithm either lowers or preserves the  average distortion. Consider the infinite sequence of average distortions obtained by repeatedly  applying the centering algorithm. Such a sequence is non-increasing, and is also bounded below by zero. Hence this sequence converges. However, it could be that  the limiting average distortion is a local minimum. 

 We show that log-concavity of the density implies that a local minimum must be a global minimum, 
 in the special case where the random variable~$X$ is scalar, and the distortion function~$ \delta\left(\cdot\right)$ is 
 the usual  square error, or the  absolute error.

 The technique of our proofs 
 is based on the generic family of all intervals with a prescribed probability mass. We span this family with a sliding interval, that has a variable length but  a fixed probability mass. We shall compare the speed at which the midpoint of this interval moves, when compared to speeds of the conditional mean and the median of the interval. 
\subsection{Centering minimizes the mean squared error}
We first need to establish an inequality relating the values of a log-concave density at end points of a finite interval, to the  average of the density over the interval. Over an interval~$ \left[ a , b  \right]$ the average of the density is defined as:
\begin{align}
    p_{\text{average}} & \triangleq 
    {\frac{1}{\left( b - a \right) }} 
    \int_{a}^{b} p_X(x) dx  .
    \label{eqn:lebesgueAverageDefinition}
\end{align}
Because of log-concavity, both the values~$ p_X(a) $ and $p_X(b)$ cannot be simultaneously  higher than~$ p_{\text{average}}  .$ In fact we get the following lower bound on the average.
\begin{lemma}%
    Consider a finite interval~$ \left[ a , b  \right]$ on which a given density~$p_X{(\cdot )}$ is log-concave and has a positive probability mass~$\eta$, with~$ p_{\rm{average}} $  defined by~\eqref{eqn:lebesgueAverageDefinition}. Then it follows that:
\begin{align}
    p_{\rm{average}}    & \ge
    \sqrt{%
           p_X \left( a \right)   \cdot   p_X \left( b \right)    
         }
 .
    \label{eqn:logarithmicInequality}
\end{align}
\end{lemma}
\begin{proof}
The logarithm function is concave, and the density function~$ p_X\left( \cdot \right)$ is non-negative. Applying the integral form of Jensen's inequality on the $\log$ function, we get:
    \begin{align*}
        \log \left(  
              {\frac{1}{\left( b - a \right) }} 
               \int_{a}^{b} p_X(x) \, dx  
            \right) 
        & \ge 
              {\frac{1}{\left( b - a \right) }} 
               \int_{a}^{b} \log 
               p_X(x)
               \, dx  .
    \end{align*}
By the concavity of the graph of~$ \log{ p_X(\cdot)}$, it follows that  between $x = a $ and $x = b , $  the area under this curve is bigger than the area under the chord that connects the point~$ \left( a , \log{p_X(a)} \right)  $ to the point~$ \left( b , \log{p_X(b)} \right)  .$ Using this and the definition of the average, the last inequality leads us to:
    \begin{align*}
        \log{ 
            p_{\text{average}} 
            } 
        & \ge 
              {\frac{1}{2}} 
                \Bigl( \log{p_X(a)} + \log{ 
                     p_X(b) 
                     } 
                \Bigr) 
                .
    \end{align*}
Exponentiating on both sides leads to Inequality~\eqref{eqn:logarithmicInequality}.
\end{proof}

\begin{theorem}
    Let~$X$ be a scalar random variable with a density~$p_X{\left( \cdot \right)}$ that is logarithmically concave. Then given any probability value~$\eta \in \left[ 0 , 1 \right] $, either:
    \begin{itemize}
        \item{there is
a unique interval~$ 
{\mathscr{S}}^* \triangleq \left[ a^*_\eta , a^*_\eta + l^*_\eta \right]
$
with probability mass of~$\eta,$ and minimizing the conditional variance~$ 
    {\mathbb{E}}\left[ 
            {\left( X  -  
                    {\mathbb{E}}\left[ X  \left\lvert   
                                       X \in S
                                          \right.
                               \right]
            \right)}^2
    \left\lvert  
            X \in S
    \right.
    \right]
          , $
over all silence sets~${\mathscr{S}}$ having probability mass at least~$\eta$, or%
}
        \item{there is a sliding family containing an infinity of optimal intervals all of the same length. In specific, there is
a unique positive length~$l_{\eta}^*,
$
a unique lower bound~$ \underline{a}_\eta ,
$ and
a unique upper bound~$ \overline{a}_\eta
$
such that, for every left end~$ a $ within these bounds, the 
interval~$ 
 {\mathscr{S}}_a^*  \triangleq  
\left[ a ,  a + l_\eta^* \right]
$
has a probability mass of~$\eta,$
and minimizes the conditional variance
over all silence sets~${\mathscr{S}}$ having
a probability mass of at least~$\eta$.%
}
    \end{itemize}
    \label{thm:meanSquareError}
\end{theorem}
\begin{proof}
    It is enough to only consider silence sets that are intervals, because the centering algorithm applied on a silence set for a log-concave density produces an interval after  only one run. And we know that this silence interval has an average distortion that is no worse than that of the silence set before the run.

    For scalar random variables, we define a density's {\textit{tail}} as a fragment, over an interval that can be finite or infinite, such that the density is (i)~non-zero over this interval (except possibly at ends of the interval), 
    (ii)~tends to a limit of zero at one end of the interval (the density is allowed to be discontinuous at this point, but at least a one-sided  limit exists  with a limiting value of zero), (iii)~reaches the density's peak at the other end, and (iv)~is monotonically increasing when we proceed from the zero-limit end to the peak end. Consider for example the exponential density~$ \lambda e^{-\lambda x} $ for nonegative~$x.$  
    This density has only one tail, which is the whole waveform of the density. 
    Except for the uniform density on a finite interval, every scalar log-concave density has at least one tail.
    The Gaussian, Cauchy or Rayleigh densities etc. have two tails.

    We separately treat three cases: (1)~the uniform density which is the only log-concave density that has no tail, (2)~densities that have tails on only one side of their peaks, and (3)~densities that have tails on both sides.

\noindent
{\textsf{\textbf{Case 1~Uniform density:}}}
For any probability mass~$\eta$ that is positive but less than one, the corresponding silence set is non-unique. Any interval is optimal, whose probability mass equals the given value of~$\eta.$

\noindent
{\textsf{\textbf{Case 2 Density with tail on only one side of peak:}}}
 \begin{figure}
   \centering
\includegraphics[width=0.4\textwidth]{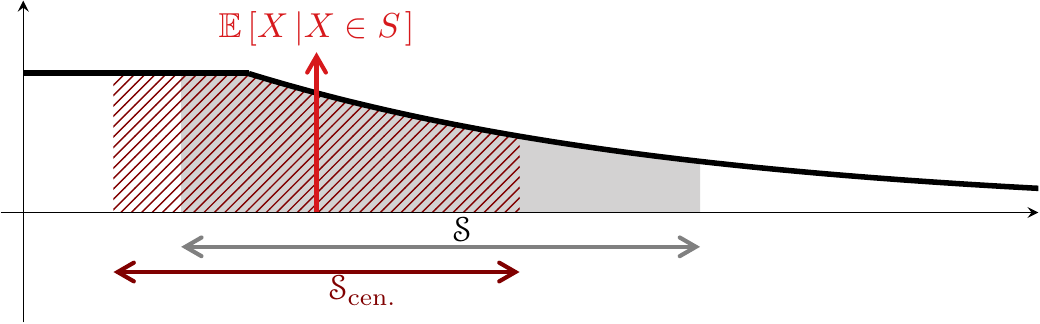}
\caption{Case~2 density with tail on only one side.}
\label{fig:oneSidedDensity}
\end{figure}
Without loss of generality, we shall assume that the tail is to the right of the peak as shown in Figure~\ref{fig:oneSidedDensity}. 

A {\textit{falling interval}} is any interval~${\mathscr{F}} = \left[ a , b \right]$ such that the density is non-increasing when we move from the left end~$a$ to the right end~$b$, and the density value at the right end is strictly less than at the left end. For example, both the intervals~${\mathscr{S}} ,  {\mathscr{S}}_{cen.} $ of Figure~\ref{fig:oneSidedDensity} are falling intervals.

On every falling interval~${\mathscr{F}}$, the probability mass of the first half of the interval namely~$\left[ a , ( a + b ) / 2 \right] ,$ is strictly greater than the probability mass of the second half~$\left[ ( a + b ) / 2 , b \right] .$ Moreover, on any sub-interval of~${\mathscr{F}}$  that is symmetric around the midpoint~$ ( a + b ) / 2 $ the density:
\begin{itemize}
    \item{either has the same probability mass on either half of the sub-interval - this happens only if the sub-interval  lies entirely under the flat top of the density, or,}
    \item{has more mass on the  left side of the sub-interval.}
\end{itemize}
Moreover we can be sure that are indeed sub-intervals with more mass on their left halves. In specific there is a non-negative lower bound~$l$ that is strictly less than~$ ( a + b ) / 2  $ such that, every  sub-interval that has length longer than~$2 l$ and is symmetric about~$ ( a + b ) / 2  $ does have more probability mass in its left half. 
Hence the conditional mean of every falling interval must be to the left of its midpoint. 

Suppose that 
a falling interval~${\mathscr{F}}$ is such that the density is non-zero in at least a small neighbourhood to the left of the left end of~${\mathscr{F}}$. Then the centering~${\mathscr{F}}_{cen.}$ must have its conditional mean strictly to the left of the conditional mean of~${\mathscr{F}}$ itself.  

Suppose on the other hand, that 
a falling interval~${\mathscr{F}}$ is such that the density is zero everywhere to the left of the left end of~${\mathscr{F}}$. Then the centering~${\mathscr{F}}_{cen.}$ 
equals~${\mathscr{F}}$ itself. 

We shall now consider a sliding interval  
of probability mass exactly~$\eta .$ As the left end 
varies, we shall consider 
the graphs of the conditional mean and the midpoint. 
We shall describe the nature of intersections between these two  graphs, 
using our above  conclusions about falling intervals.

If no flat top exists for a one-tailed density 
with the tail to the right of the peak, 
then the graphs of the conditional mean and the midpoint can have no intersection. 
The only 
possible graphical behaviour is:
(ii-a)~
the graph of the conditional mean starts below the graph of the midpoint, and as the left end~$a$ increases, the gap
between the conditional mean and  the midpoint keeps  increasing.

If a flat top does exist for a one-tailed density with the tail to the right of the flat top, then the graphs of the conditional mean and the midpoint could  have  either no intersection or coincide over an interval. There are two possible graphical behaviours. First is the behaviour~(ii-a) listed in the previous paragraph, and this happens only if the probability mass~$\eta$ is large enough that the only intervals reaching this mass are falling ones. The other possibility is: (ii-b)~in which the conditional mean coincides with the midpoint over a finite interval~(which is a sub-interval of the support of the flat top), and to the right of this interval, the conditional mean is lesser than the midpoint - this possibility happens when the probability mass~$\eta$ is small enough so that there are corresponding intervals that fall entirely within the support of the density's flat top.
\begin{figure}
    \centering
    \subfloat[Exponential~($\lambda e^{-\lambda x},$ for $x \ge 0$) with $\lambda = 1,$ $\eta  = 0.75$]{%
\includegraphics[width=0.45\textwidth]{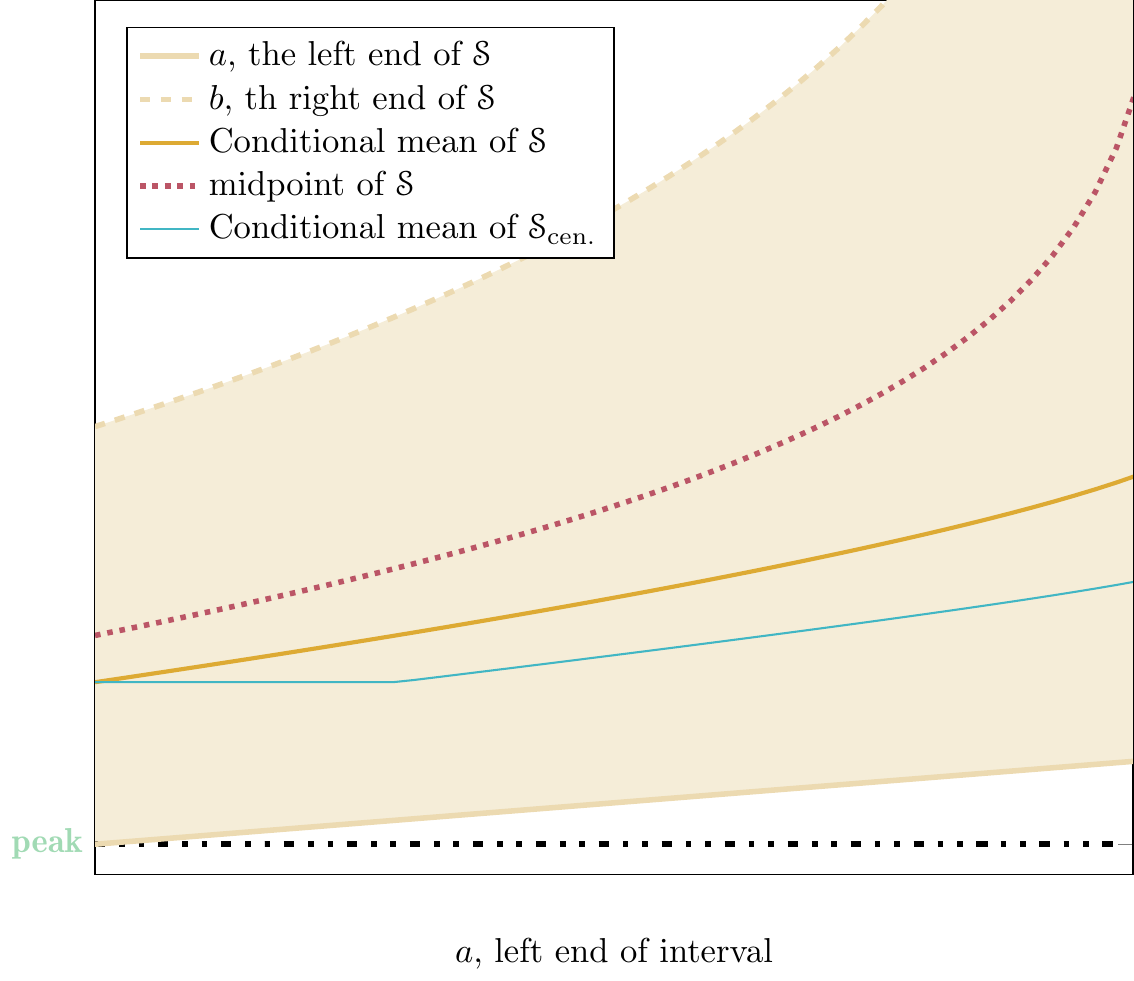}
\label{subfig:conditionalMeansMeetAtLeftEnd}
}
\linebreak
    \subfloat[Rayleigh~($( x / \sigma^2 ) e^{ - {\tfrac{x^2}{2\sigma^2}}},$ for $x \ge 0 )$ with~$\sigma = 8$, and  $ \eta  = 0.8$]{%
\includegraphics[width=0.45\textwidth]{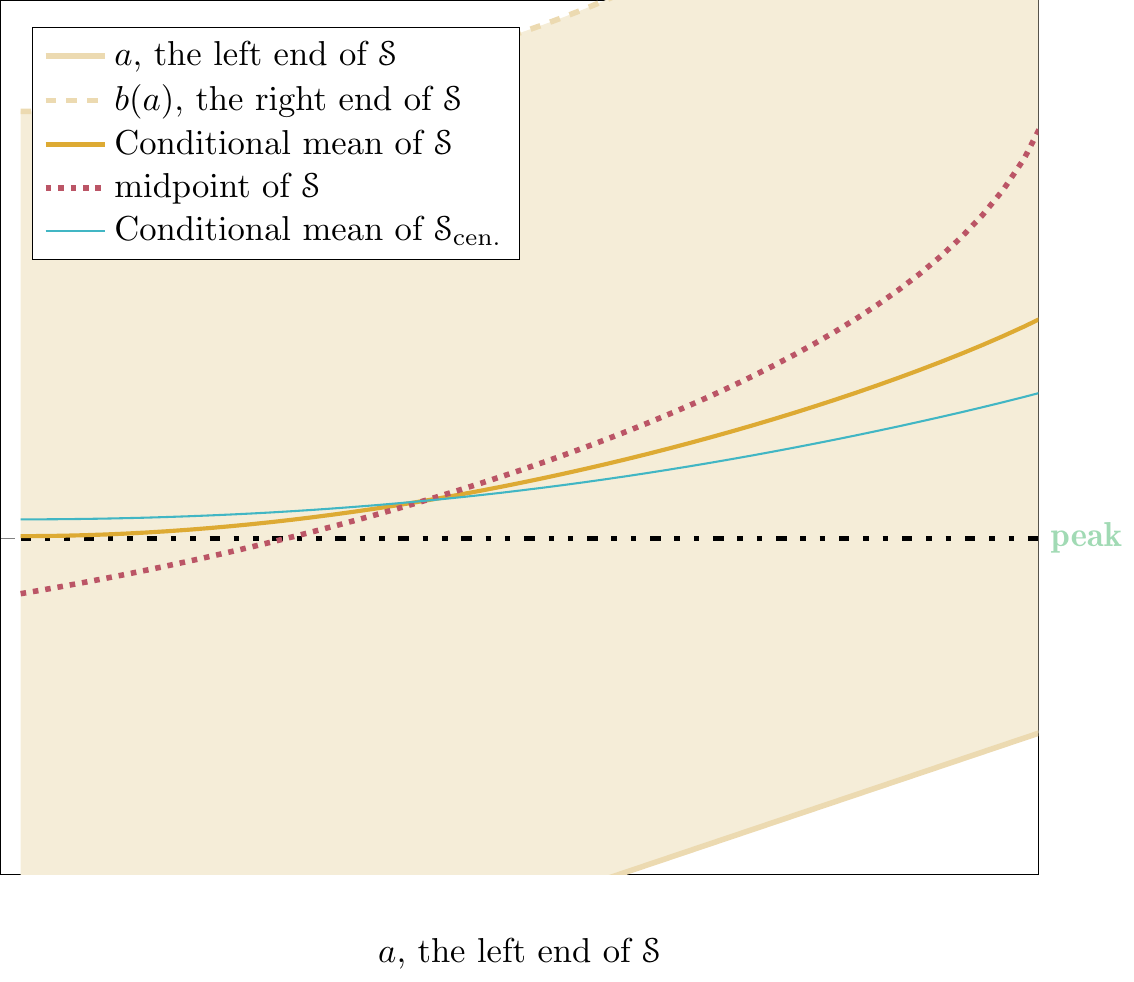}
\label{subfig:conditionalMeansMeetNearPeak}
}
\linebreak
\subfloat[Rayleigh with an artificially pasted flat top, and  $ \eta = 0.1$]{%
\includegraphics[width=0.45\textwidth]{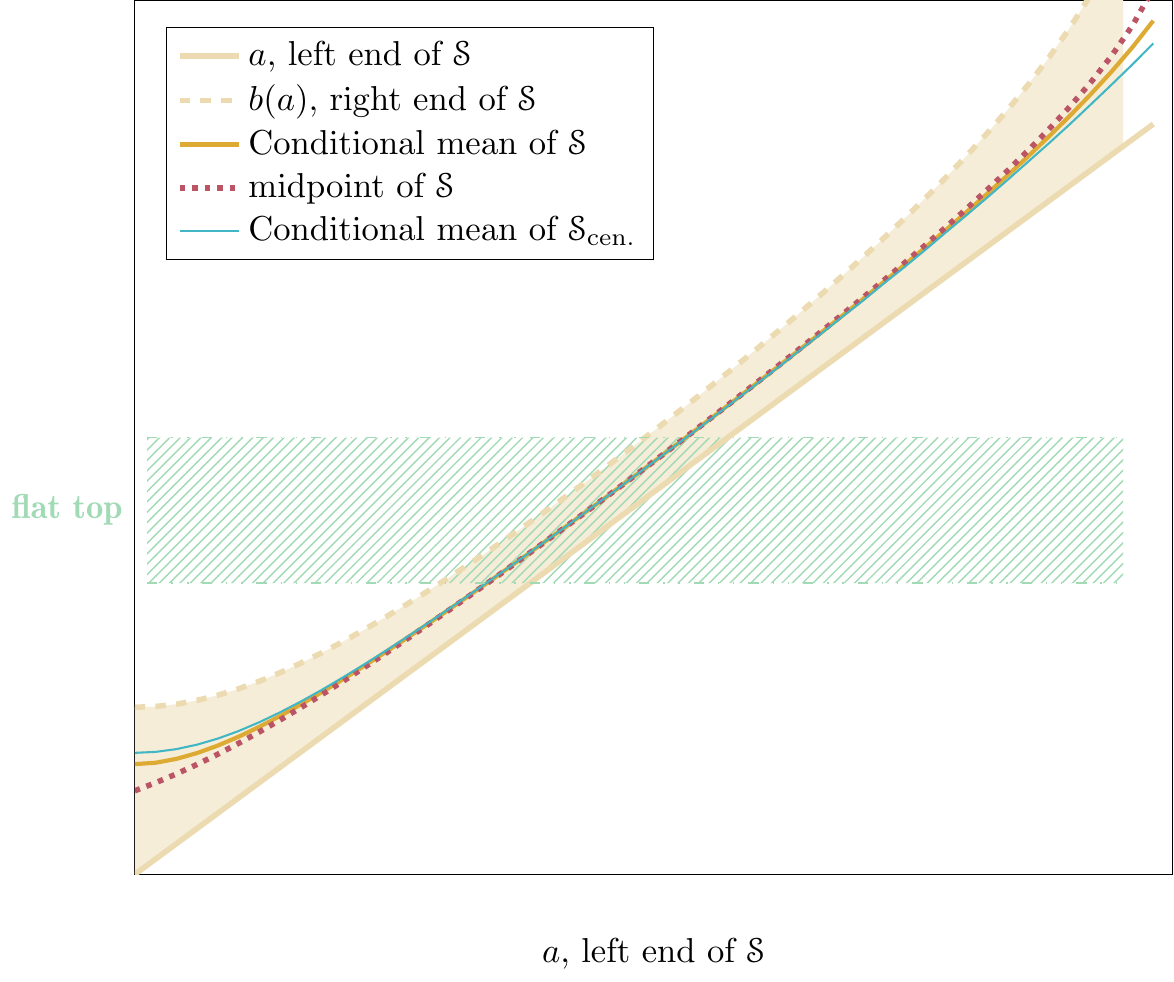}
\label{subfig:conditionalMeansCoincideOverInterval}
}
\caption{The different prototypical ways in which the conditional means of~${\mathscr{S}}$ and of~${\mathscr{S}}_{\text{cen.}}$ can intersect}
\label{fig:howConditionalMeansMeet}
\end{figure}
Only two kinds of centered intervals are possible 
for a density that has a tail only on the right hand side: 
a falling interval that starts at the left end of the support of the density, or
an interval where the conditional mean equals the midpoint. 
Under possibility~(ii-a) only the first kind of centered interval exists and is unique, and under possibility~(ii-b) only the second kind exists and can be non unique.

\noindent
{\textsf{\textbf{Case 3 Density with tails on both sides of peak:}}}
In handling Case~2 we had defined and studied the notion of falling intervals. Similarly we shall define a {\textit{rising interval}} as 
any interval~${\mathscr{I}} = \left[ a , b \right]$ such that the density is non-decreasing when we move from the left end~$a$ to the right~$b$ end, and the density value at the right end is strictly greater than at the left end. A {\textit{peaked interval}} is any having  an interior point where the density is greater than at either ends of the interval.

For a two-tailed density, we shall show that an optimal interval can neither be a rising interval nor a falling one. Rather it has to be either a peaked interval or an interval over which the density is uniform.

We shall study the sliding family of silence intervals~${\mathscr{S}}_a = [a, b(a)] $ parametrized by the variable left end~$a,$ and possessing the property that the probability mass is~$\eta .$ The variable~$b$ is implicitly defined via the probability mass constraint:
\begin{align*}
    \int_{a}^{b} p_X(x) dx & = \eta .
\end{align*}
When the LHS above is viewed as a function of two independent variables~$a, b$ the directional derivative vector is:~$ \begin{bmatrix} p(a) & p(b) \end{bmatrix} .$ If the left end~$a$ lies in the interior of the support of the density, then the derivative of~$b$ with respect to~$a$ exists. By the implicit function theorem, under this condition, the right end~$b(a)$ must be a continuously differentiable function of the left end~$a.$
And therefore the conditional mean and midpoint of~${\mathscr{S}}_a$ vary as continuously differentiable functions of~$a,$ if~$p_X\left(a\right)$ is positive. 

For what values of the left end~$a$ does the conditional mean of~${\mathscr{S}}$ equal its midpoint~? Note that for such values of~$a$, the silence interval~${\mathscr{S}}$  
equals its centering.  We shall now show that for a log-concave density, the set of such values of~$a$ can only be an empty set, or a singleton set, or an finite interval.

What are the speeds at which the conditional mean and the midpoint of~${\mathscr{S}}$ increase as we increase~$a$? The ratio is:
\begin{gather*}
    {\frac{
            {\frac{1}{2}} 
            {\frac{d}{da}} \left( {a} + {b(a)}  \right)
          }
          {
            {\frac{1}{\eta}} 
            {\frac{d}{da}} \int_{a}^{b(a)}{ x p_X (x) dx} 
          }
    }
    \ = \
    {\frac{\eta}{2}} 
    {\frac{ p_X\left(a\right) + p_X\left( b \right) }
          { \left( b- a \right) p_X\left(a\right)  p_X\left( b \right) }
    } ,
\\
    \ \ge \
    {\frac{\eta}
          { \left( b- a \right) \sqrt{ p_X\left(a\right)  p_X\left( b \right) } }
    }
    \quad (\text{by the AM-GM inequality}) ,
\\
    \ = \
    {\frac{ p_{\text{average}}  }
          { \sqrt{ p_X\left(a\right)  p_X\left( b \right) } }
    } ,
    \phantom{abcdefghijklmnoabcdeaafghijklmno}
\\
    \ \ge \ 1
     \quad \quad (\text{by log-concavity~\eqref{eqn:logarithmicInequality}})
 .   \phantom{abcdefghijklmnnoabcbcz}
\end{gather*}
In fact the ratio on the LHS above can equal the upper bound of one, only if~$ p_X\left(a\right) = p_{\text{average}} =  p_X\left( b \right) .$ And because of log-concavity, this double equality is the same as the density being uniform over the interval~$ \left[ a , b \right] . $
In summary, as we increase the left end~$a,$:
\begin{itemize}
    \item{either the conditional mean increases at a rate that is strictly slower than that of the midpoint, or,}
    \item{the conditional mean increases at exactly the same rate as the midpoint, and this can only happen when the density has a flat top and the entire interval~$ \left[ a , b \right] $ is within the support of this flat top. When this happens, clearly the conditional mean is the same as the midpoint.}
\end{itemize}
Recall that in a falling interval  the conditional mean is to the left of the midpoint, and that in a rising interval  the conditional mean is to the right of the midpoint. In a peaked interval, the conditional mean can be on either side of the midpoint, or equal it.

For a two-tailed density, if the silence probability~$\eta$ is close enough to zero, then there do exist corresponding  silence intervals of every kind: falling, rising and peaked. But if the silence probability~$\eta$ is close enough to one, then it could happen that no corresponding  silence intervals  exists of the falling kind, or of the rising kind, or of neither of these two kinds.
From these considerations, we shall systematically list the  nature of possible intersections of the graphs of the conditional mean and of the midpoint. 

If no flat top exists for a two-tailed density, then it has a unique peak and so the graphs of the conditional mean and the midpoint can have at most one intersection, and moreover such an intersection has to be transversal (no tangency). Thus the following three are the only possible crossing behaviours of the conditional mean and the midpoint, as the left end~$a$ varies:
(iii-a)~
the graph of the conditional mean starts above the graph of the midpoint, crosses it once (transversally), and ends up below it.
(iii-b)~
the graph of the conditional mean starts below the graph of the midpoint and simply keeps on separating more from it,
(iii-c)~the graph of the conditional mean starts above the graph of the midpoint and keeps moving towards it but never manages to reach it.

If a flat top does exist, then the graphs of the conditional mean and the midpoint can either have no intersection, or have exactly one intersection, or coincide  over an interval. This allows four  possible crossing behaviours of the conditional mean and the midpoint. 
The first three are precisely those listed in the previous paragraph, namely: (iii-a), (iii-b), and (iii-c). The fourth possibility is:
(iii-d)~
the graph of the conditional mean starts above the graph of the midpoint, meets it and coincides with it over a finite interval~(which lies strictly within the interval supporting the flat top), and to the right of this interval the conditional mean lies  below the midpoint. 

Three kinds of centered intervals are possible 
for a two-tailed density:
a peaked interval that starts at the left end of the support of the density, or
an interval where the conditional mean equals the midpoint, or 
a peaked interval that ends at the right end of the support of the density.
Under possibility~(iii-a) and~(iii-d) only the first kind of centered interval exists -  unique for possibility~(iii-a), and likely non unique for possibility~(iii-d). Under~(iii-b) only the first kind exists and is unique, and under~(iii-c) only the third kind exists and is unique.
\end{proof}
\begin{corollary}
A centered silence set minimizes the mean squared estimation error, if the density of the sampled random variable is logarithmically concave.
\label{corollary:centeringMinimizesMeanSquaredError}
\end{corollary}
\begin{proof}
Any non-increasing, and positive sequence of numbers must converge.
Therefore the process of repeated iteration of the centering operation must lead to a limiting value for the conditional distortion.

An iteration of the centering operation cannot increase the distortion. 
If such an iteration results in the same distortion as before, then the original and newer sets can differ only by a set of measure zero. In which case the newer set must be a centered one.
The Corollary follows from the guarantee of Theorem~\ref{thm:meanSquareError}
that centered sets are essentially unique.
\end{proof}
\subsection{Centering minimizes the mean absolute error}
Given the random variable~$X,$ its density~$p_X\left(\cdot\right),$ and the interval~$  [  a , b ]$, let ~$\mu_{[a,b]}$ denote the median of the random variable~$X $ conditioned on this interval. That is to say that:
\begin{align}
\int_a^{ \mu_{[a,b]} }
{ p_X \left( x \right)  dx }
=
\int^b_{ \mu_{[a,b]} }
{ p_X \left( x \right)  dx } 
\  =  \  
{ \frac{1}{2}} \cdot 
{\mathbb{P}} \left[ X \in \left[ a , b \right] \right]
.
\end{align}
Next we prove an inequality for log-concave densities,   similar to Lemma~\ref{eqn:logarithmicInequality}. We shall use it when the distortion measure is the absolute value of the error.
\begin{lemma}%
    Consider a finite interval~$ \left[ a , b  \right]$ on which a given density~$p_X{(\cdot )}$ is log-concave.
    Then the value of the density~ at the conditional median~$  
      \mu_{[a,b]}
    $
    is no less than the square root of the values of the density at the end points:~$ a , b .$ In other words:
\begin{align}
    p_X \left(
        { \mu_{[a,b]} } 
        \right)
        & \ge
    \sqrt{ 
           p_X \left( a \right) \cdot    p_X \left( b \right)    
         }    .
\label{eqn:medianInequality}
\end{align}
\label{thm:mediamInequality}
\end{lemma}
\begin{proof}
Consider the graph of the logarithm of the given density function, over the interval~$[ a , b ]$. As show in Figure~\ref{fig:densityAtMedianHigherThanGeometricMeanOfEnds}, 
 over the interval~$[ a , b ] , $ log-concavity 
implies that 
the graph of the logarithm of the density cannot go below the  corresponding chord.
\begin{figure}
\begin{center}
\begin{tikzpicture}[scale=0.9]
\begin{axis}[
        xtick={0.5,1.65,2.32,2.8},
        xticklabels={$a$,${{\left(a + b\right)}/2}$,$\nu$,$b$},
        ytick={0.76,2.255,3.75},
        yticklabels={$\log{p\left(b\right)}$,%
                     $\log{\sqrt{ p\left(a\right) p\left(b\right) }}$,%
                     $\log{p\left(a\right)}$,}%
    ]
    \addplot [draw=none, name path = xaxis, domain=0:2.91,]  {0};
    \addplot [ draw = none, domain=0.2:2.9, samples=12, smooth, name path = logOfDensity,]
    {4-(x-1)^2};
    \addplot [ 
        draw=none,
        pattern color=blue2,
        fill=blue2,
        blue2,
        pattern= north east lines, 
    ]
    fill between[
        of = logOfDensity and xaxis,
        soft clip={domain=0.5:2.32},
    ];
        \addplot [ultra thick, color=black, domain=0.2:3.1, samples=12, smooth, name path = logOfDensity,]
    {4-(x-1)^2};
    \addplot [ 
        draw=none,
        pattern color=brown2,
        fill=brown2,
        brown2,
        pattern= north east lines, 
    ]
    fill between[
        of = logOfDensity and xaxis,
        soft clip={domain=2.32:2.8},
    ];
    %
        \addplot [ultra thick, color=black, domain=0.2:2.9, samples=12, smooth, ]
    {4-(x-1)^2};
\addplot[thick, dotted ,color=blue5, domain=0:0.5, samples=3, smooth,]
    {3.75};
    \addplot[thick, dotted ,color=blue5, domain=0:2.8, samples=5, smooth,]
    {0.76};
    \addplot[very thick, dashed ,color=red1, domain=0:2.32, samples=5, smooth,] {2.255};
    \addplot +[mark=none, color=red1, dashed] coordinates {(1.65, 0) (1.65, 2.255)};
    \addplot +[mark=none, color=red1, dashed] coordinates {(2.32, 0) (2.32, 2.255)};
\addplot[thick, color=gray, domain=0.5:2.8, samples=5, smooth, name path = chord]
    {3.75-1.3*(x-0.5)};
    \node[] at (0.5,3.75) {${\Large{\mathbf{\bullet}}}$}; 
    \node[] at (2.8,0.76) {${\Large{\mathbf{\bullet}}}$}; 
    \node[color=red1] at (1.65,2.255) {${\large{\mathbf{\bullet}}}$}; 
    \node[color=red1] at (2.32,2.255) {${\large{\mathbf{\bullet}}}$}; 
\end{axis}
\end{tikzpicture}
\end{center}
\caption{Illustrating that the density at the median is not lower than the geometric mean of the density values at the end points%
\label{fig:densityAtMedianHigherThanGeometricMeanOfEnds}}
\end{figure}
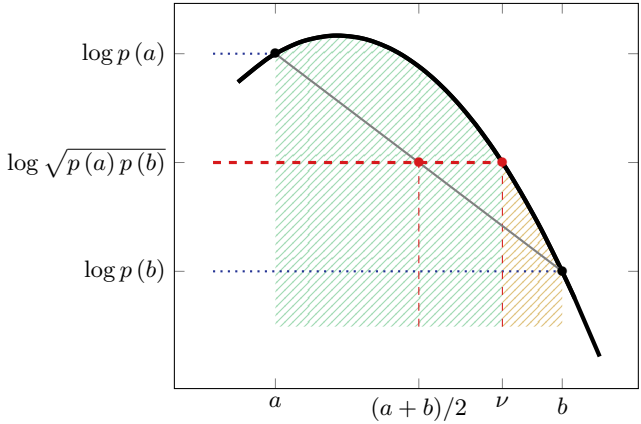

We have the following three mutually exclusive cases:

{\textbf{Case~i:~$ p_X (a) = p_X (b) $.}}
Log-concavity implies unimodality. Hence for any real~$x $ in the open interval~$ \left(    a , b \right) , $ the density at that point 
cannot be lesser than the minimum of the values at the ends of the interval. Hence the Lemma holds in this case.

{\textbf{Case~ii:~$ p_X (a) > p_X (b) $.}}
Consider the point~$ \nu $ from within the interval~$ [ a , b ] $ such that:
\begin{align*}
    p_X \left(  \nu \right)
    & =
    \sqrt{%
     p_X \left(  a \right)
      \cdot p_X \left(  b \right)
    }
\end{align*}
As shown in Figure~\ref{fig:densityAtMedianHigherThanGeometricMeanOfEnds}, we get the
 following geoemetric facts:
\begin{itemize}
    \item since~$ p_X \left(  a \right)$ is greater than~$ p_X \left(  b \right)$ it follows that
    \begin{gather*}
          \left(  a + b \right) / 2 
          \  \le \ 
          \nu  
          \ \le \ 
          b ,
    \end{gather*}
    In other words, the length of the interval~$   \left[
          \nu  , b \right]
    $
    cannot be larger than that of the interval~$   \left[
         a  ,
          \nu \right]
      . $
    \item 
      the value of the density~$ p_X \left(  \cdot  \right) $ at every point in the interval~$   \left[
          \nu  , b \right]
    $
    is lesser or equal to the value of the density at every point of the interval~$   \left[
         a  ,
          \nu \right]
      . $
\end{itemize}
Hence the integrals over these two intervals must obey:
\begin{align*}
\int_{\nu}^{b}{ p_X  \left( x \right) dx }
    &
    \le
  \int_{a}^{\nu}{ p_X  \left( x \right) dx }  
  .
\end{align*}
Hence the median~$  \mu_{ [ a , b ] }$ cannot lie to the right of the point~$ \nu . $ Hence the Lemma follows in this case also.

{\textbf{Case~iii:~$ p_X (a) < p_X (b) $.}}
This case can be handled exactly as in Case~ii, except that the point~$\nu$ lies to the left of~$  \left(  a + b \right) / 2 .$
\end{proof}

\begin{theorem}
    Let~$X$ be a scalar random variable with a density~$p_X{\left( \cdot \right)}$ that is logarithmically concave. Then given any probability value~$\eta \in \left[ 0 , 1 \right] $, either:
    \begin{itemize}
        \item{there is
a unique interval~$ 
{\mathscr{S}}^* \triangleq \left[ a^*_\eta , a^*_\eta + l^*_\eta \right]
$
with probability mass of~$\eta,$ and minimizing the conditional absolute error~$ 
    {\mathbb{E}}\left[ 
       \left.
            {\left\lvert X  -  mu_{  [ a , b  ] }
            \right\rvert}
    \right\rvert  
            X \in S
    \right]
          , $
over all silence sets~${\mathscr{S}}$ having
a probability mass of at least~$\eta$, or%
}
        \item{there is a sliding family containing an infinity of optimal intervals all of the same length. In specific, there is
a unique positive length~$l_{\eta}^*,
$
a unique lower bound~$ \underline{a}_\eta ,
$ and
a unique upper bound~$ \overline{a}_\eta
$
such that, for every left end~$ a $ within these bounds, the 
interval~$ 
 {\mathscr{S}}_a^*  \triangleq  
\left[ a ,  a + l_\eta^* \right]
$
has a probability mass of~$\eta,$
and minimizes the conditional absolute error
over all silence sets~${\mathscr{S}}$ having
a probability mass of at least~$\eta$.%
}
    \end{itemize}
    \label{thm:meanAbsoluteError}
\end{theorem}
\begin{proof}
The proof follows almost exactly like that for Theorem~\ref{thm:meanSquareError}. The only difference is that instead of  the silence interval's conditional mean, we use the median~$ \mu_{  [ a, b ]  } . $

By differentiating both sides of the 
 chance constraint:
\begin{align*}
    \int_{a}^{ \mu_{  [ a, b ]  }  } p_X\left( x \right) dx 
 & =
 { \frac{1}{2}} \eta   
 \\
 \ \text{we get:}  \ \    \ \ 
{ \frac{d}{da}}  \mu_{  [ a, b ] }
& ={ \frac{
             p_X\left( a \right) 
          }
          { 
             p_X\left(  \mu_{  [ a, b ] } \right) 
          }} .
\end{align*}

What are the speeds at which the median and the midpoint of~${\mathscr{S}}$ increase as we increase~$a$~? The ratio is:
\begin{gather*}
    {\frac{
            {\frac{1}{2}} 
            {\frac{d}{da}} \left( {a} + {b(a)}  \right)
          }
          {
            {\frac{1}{\eta}} 
            {\frac{d}{da}} \mu_{  [ a, b ] } 
          }
    }
    \ = \
    {\frac
          { 
            \left( p_X\left(a\right) + p_X\left( b \right) \right)
             p_X\left( \mu_{  [ a, b ] 
          } \right) 
          }
          { 2 p_X\left(a\right)  p_X\left( b \right) 
    } },
\\
    \ \ge \
    {\frac{    p_X\left( \mu_{  [ a, b ]  } \right) }
          { \sqrt{ p_X\left(a\right)  p_X\left( b \right) } }
    }
    \quad (\text{by the AM-GM inequality}) ,
\\
    \ \ge \ 1 ,
     \quad \quad (\text{by log-concavity~\eqref{eqn:medianInequality}})
 .   \phantom{abcdefghijkoabcbcz}
\end{gather*}
The rest of the proof is the same as for Theorem~\ref{thm:meanSquareError}.
\end{proof}
\begin{corollary}
A centered silence set minimizes the mean absolute estimation error, if the density of the sampled random variable is logarithmically concave.
\end{corollary}
\noindent
The proof is similar to that for Corollary~\ref{corollary:centeringMinimizesMeanSquaredError}.
\section{Super-level intervals are nearly optimal\label{section:superLevelInterval}}
We shall calculate the conditional variance of the following families of silence intervals. These famlilies result from heuristic attempts to keep the interval lengths as small as possible, while collecting the required probability mass~$\eta$:
\begin{itemize}
    \item{{\textbf{Super~level~sets:}} A super level interval is defined as:
        \begin{align*}
            {\mathscr{S}}^{\text{super-level}}_{\eta} & = 
                      \left\{
                              x \in {\mathbb{R}}^n 
                         : 
                             p_X \left( x - \mu 
                             \right) \ge \alpha_{\eta}^{\text{super-level}}
                          \right\},
        \end{align*}
        where the level~$\alpha^{super-level}_{\eta}$ is chosen to be the smallest level guaranteeing that the above interval has a probability mass of at least~$\eta.$
        }
    \item{{\textbf{Equal sides around mode:}} An interval with equal sides around the mean is defined as: 
        \begin{align*}
            {\mathscr{S}}^{\text{equal-sides}}_{\eta} & = 
                      \left\{
                              x \in {\mathbb{R}}^n 
                         : 
                             \left\lvert  
                             x - \mu 
                             \right\rvert  
                             \le \alpha_{\eta}^{\text{equal-sides}}
                          \right\},
        \end{align*}
        where the half interval width~$\alpha^{equal-sides}_{\eta}$ is chosen to be the smallest one guaranteeing that the above interval has a probability mass of at least~$\eta.$
        }
    \item{{\textbf{Equal areas around mode:}} An interval with equal areas around the mean is defined as:~$
   {\mathscr{S}}^{\text{equal-areas}}_{\eta}  =
     \left[  \alpha_{\eta}^{\text{equal-areas}}, \beta_{\eta}^{\text{equal-areass}}\right], $ where the limits are such that
        \begin{gather*}
\int_{ \alpha_{\eta}^{\text{equal-areas}}}^{\mu}{  p_X (x) dx}  \ =  \
\int_{ \mu }^{\beta_{\eta}^{\text{equal-areas}}}{  p_X (x) dx}  \ =  \
            {\frac{\eta}{2}}  .
        \end{gather*}%
}
    \item{{\textbf{Mode-as-conditional~mean:}} An interval with the conditional mean as its mode is defined as:~$
   {\mathscr{S}}^{\text{mode-as-mean}}_{\eta}  =
     \left[  \alpha_{\eta}^{\text{mode-as-mean}}, \beta_{\eta}^{\text{mode-as-mean}}\right], $ 
where the limits are carefully chosen around the mean  such that
        \begin{align*}
                    {\mathbb{E}}\left[ X  \left\lvert   
                                       X \in 
     \left[  \alpha_{\eta}^{\text{mode-as-mean}}, \beta_{\eta}^{\text{mode-as-mean}}\right] 
                                          \right.
                               \right]
                               & = \mu .
        \end{align*}
        }
\end{itemize}

 \begin{figure}
\centering
\subfloat[Laplace modified to have unbalanced sides]{%
 \begin{tikzpicture}
     \begin{axis}[
             width = \linewidth, height = 0.4\linewidth,
          xmin=-4.5, xmax=4.5,
          ymin=-0.05, ymax=1.35,
          axis lines=center,
          axis on top=true,
          xlabel={},
          ylabel={},
          ticks = none,
          axis on top,
      ]
     \addplot [name path =pdfCurveLeft, mark=none,draw=black,ultra thick, samples = 200, domain=-4.5:0,] {exp(1.5*\x)};
     \addplot [name path =pdfCurveRight, mark=none,draw=black,ultra thick, samples = 200, domain=0:4.5,] {exp(-0.3*\x)};
     \path[name path=xaxis] (-4.5,0) -- (4.5,0);
    \addplot [
        draw=none,
        fill=blackOne!20, 
    ]
    fill between[
        of = pdfCurveLeft and xaxis,
        soft clip={domain=-4.5:0},
    ];
    \addplot [
        draw=none,
        fill=blackOne!20, 
    ]
    fill between[
        of = pdfCurveRight and xaxis,
        soft clip={domain=0:4.5},
    ];
\end{axis}
\end{tikzpicture}
}
\linebreak
\subfloat[Patchwork of circular arcs]{%
 \begin{tikzpicture}
     \begin{axis}[
             width = 0.6\linewidth, height = 0.4\linewidth,
          xmin=-2.1, xmax=2.1,
          ymin=-0.05, ymax=1.85,
          axis lines=center,
          axis on top=true,
          xlabel={},
          ylabel={},
          ticks = none,
          axis on top,
      ]
     \addplot [name path =pdfCurveLeft, mark=none,draw=black,ultra thick, samples = 200, domain=-2:0,] {0.8*pow( 4 - \x*\x, 0.5)};
     \addplot [name path =pdfCurveRight, mark=none,draw=black,ultra thick, samples = 200, domain=0:1,] { 0.8 + 0.8*pow( 1 - \x*\x, 0.5)};
     \path[name path=xaxis] (-4.5,0) -- (4.5,0);
    \addplot [
        draw=none,
        fill=blackOne!20, 
    ]
    fill between[
        of = pdfCurveLeft and xaxis,
        soft clip={domain=-1.95:0},
    ];
    \addplot [
        draw=none,
        fill=blackOne!20, 
    ]
    fill between[
        of = pdfCurveRight and xaxis,
        soft clip={domain=0:1},
    ];
\end{axis}
\end{tikzpicture}
}
\hskip 2mm
\subfloat[Triangular density]{%
 \begin{tikzpicture}
     \begin{axis}[
          width = 0.6\linewidth, height = 0.4\linewidth,
          xmin=-1.2, xmax=1.2,
          ymin=-0.05, ymax=1.1,
          axis lines=center,
          axis on top=true,
          xlabel={},
          ylabel={},
          ticks = none,
          axis on top,
      ]
     \addplot [name path =pdfCurveLeft, mark=none,draw=black,ultra thick, samples = 200, domain=-0.25:0,] {1 + 4*  \x};
     \addplot [name path =pdfCurveRight, mark=none,draw=black,ultra thick, samples = 200, domain=0:1,] { 1 -  \x};
     \path[name path=xaxis] (-4.5,0) -- (4.5,0);
    \addplot [
        draw=none,
        fill=blackOne!20, 
    ]
    fill between[
        of = pdfCurveLeft and xaxis,
        soft clip={domain=-0.25:0},
    ];
    \addplot [
        draw=none,
        fill=blackOne!20, 
    ]
    fill between[
        of = pdfCurveRight and xaxis,
        soft clip={domain=0:1},
    ];
\end{axis}
\end{tikzpicture}
}
\caption{Log-concave densities used in our empirical study.}
\label{fig:someAsymmetricDensities}
\end{figure}
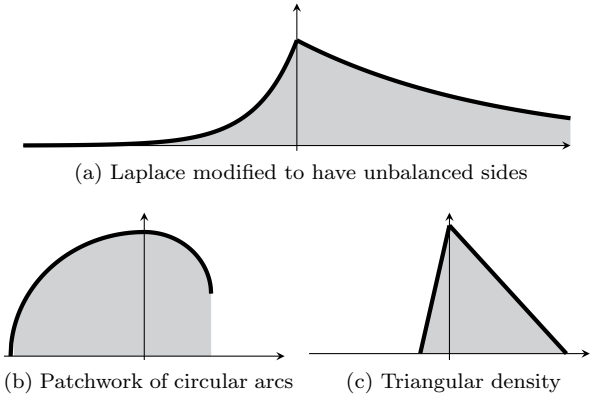
The conditional variances of these families of intervals were calculated for three log-concave densities shown in Figure~\ref{fig:someAsymmetricDensities}.  These are the unbalanced Laplace density, a density patched up from two circular arcs, and a triangular density. We took these three to be representative of log-concave density classes induced by the concavity or convexity of the waveform pieces making up the densities. 

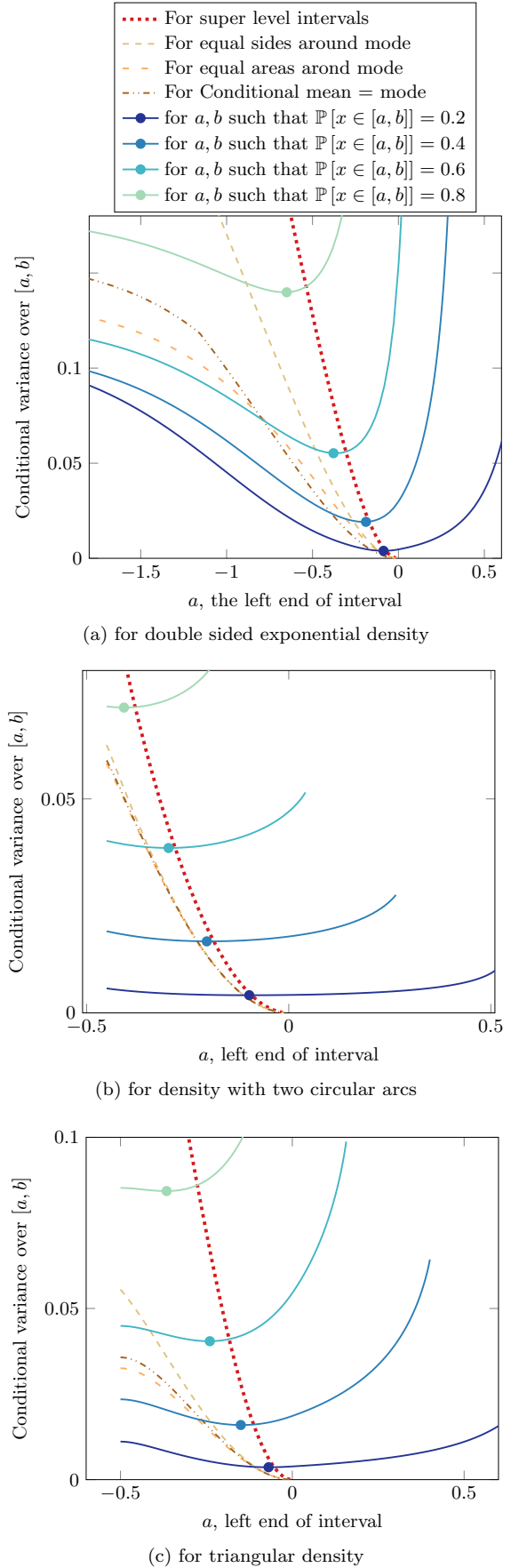
\begin{figure}
    \centering
\subfloat[for double sided exponential density]{%
\begin{tikzpicture}[scale=0.9]
\begin{axis}[ legend cell align=left,
    xlabel={${a}$, the left end of interval},
      ylabel=  {Conditional variance over~$\left[ a , b \right]$}
             ,
    xmin = -1.8, xmax= 0.6, 
    ymin = 0, ymax = 0.18, 
     xtick ={-1.5, -1, -0.5, 0, 0.5},
     ytick ={0, 0.05, 0.1, 0.15},
     yticklabels ={0, 0.05, 0.1},
    legend style ={at={(0.5,1.01)},anchor = south},
    ]    
\addplot[red1, dotted, id= forLevelSets,
 ultra thick] file {figuresAndData/dataForConditionalVarianceLaplacelevelSets.tex};
\addplot[brown2, dashed, id= forLevelSets,
 thick] file {figuresAndData/dataForConditionalVarianceLaplaceEqualSides.tex};
\addplot[red2, loosely dashed, id= forLevelSets,
 thick] file {figuresAndData/dataForConditionalVarianceLaplaceEqualAreas.tex};
\addplot[brown1, dashdotdotted, id= forLevelSets,
 thick] file {figuresAndData/dataForConditionalVarianceLaplaceMeanZero.tex};
\addplot[ mark min,
        blue5, id= forProb0pt2, thick]
   file {figuresAndData/dataForConditionalVarianceLaplaceProb0pt2.tex};
\addplot[ mark min,
        blue4, id= forProb0pt4, thick]
   file {figuresAndData/dataForConditionalVarianceLaplaceProb0pt4.tex};
\addplot[ mark min,
        blue3, id= forProb0pt6, thick]
   file {figuresAndData/dataForConditionalVarianceLaplaceProb0pt6.tex};
\addplot[ mark min,  blue2,
         id= forProb0pt8, thick]
   file {figuresAndData/dataForConditionalVarianceLaplaceProb0pt8.tex};
 \legend{ {For super level intervals},
          {For equal sides around mode},
          {For equal areas arond mode},
          {For Conditional mean = mode},
          {for $a, b$ such that ${\mathbb{P}}\left[ x \in [a , b] \right] = 0.2$},
          {for $a, b$ such that ${\mathbb{P}}\left[ x \in [a , b] \right] = 0.4$},
          {for $a, b$ such that ${\mathbb{P}}\left[ x \in [a , b] \right] = 0.6$},
          {for $a, b$ such that ${\mathbb{P}}\left[ x \in [a , b] \right] = 0.8$} };
\end{axis}
\end{tikzpicture}%
}
\linebreak
\subfloat[for density with two circular arcs]{%
\begin{tikzpicture}[scale=0.9]
\begin{axis}[ 
    xlabel={${a}$, left end of interval},
      ylabel={ Conditional variance over~$\left[ a , b \right]$  
             },
    xmin = -0.51, xmax= 0.51, 
    ymin = 0, ymax = 0.08, 
    xtick = {-0.5,  0, 0.5},
    ytick = { 0, 0.05},
    scaled ticks=false,
    yticklabels = { 0, 0.05},
    ]    
\addplot[red1, dotted, id= forLevelSets, forget plot,
 ultra thick] file {figuresAndData/dataForConditionalVarianceCircularLevelSets.tex};
\addplot[brown2, dashed, id= forLevelSets,forget plot,
 thick] file {figuresAndData/dataForConditionalVarianceCircularEqualSides.tex};
\addplot[red2, loosely dashed, id= forLevelSets,forget plot,
 thick] file {figuresAndData/dataForConditionalVarianceCircularEqualAreas.tex};
\addplot[brown1, dashdotdotted, id= forLevelSets,forget plot,
 thick] file {figuresAndData/dataForConditionalVarianceCircularMeanZero.tex};
\addplot[ mark min,  forget plot,
        blue5, id= forProb0pt2, thick] 
   file {figuresAndData/dataForConditionalVarianceCircularProb0pt2.tex};
\addplot[ mark min,  forget plot,
        blue4, id= forProb0pt4, thick] 
   file {figuresAndData/dataForConditionalVarianceCircularProb0pt4.tex};
\addplot[ mark min,  forget plot,
        blue3, id= forProb0pt6, thick] 
   file {figuresAndData/dataForConditionalVarianceCircularProb0pt6.tex};
\addplot[ mark min,  blue2,forget plot,
         id= forProb0pt8, thick ] 
   file {figuresAndData/dataForConditionalVarianceCircularProb0pt8.tex};
\end{axis}
\end{tikzpicture}%
}
\linebreak
\subfloat[for triangular density]{%
\begin{tikzpicture}[scale=0.9]
\begin{axis}[ legend cell align=left,
    xlabel={${a}$, left end of interval},
      ylabel= {Conditional variance over~$\left[ a , b \right]$ }
             ,
    xmin = -0.6, xmax= 0.6, 
    ymin = 0, ymax = 0.1, 
     xtick ={-0.5, 0, 0.5},
     ytick ={0, 0.05, 0.1},
     yticklabels ={0, 0.05, 0.1},
    ]    
\addplot[red1, dotted, id= forLevelSets,
 ultra thick] file {figuresAndData/dataForConditionalVarianceA0pt5B1pt5levelSets.tex};
\addplot[brown2, dashed, id= forLevelSets,
 thick] file {figuresAndData/dataForConditionalVarianceA0pt5B1pt5EqualSides.tex};
\addplot[red2, loosely dashed, id= forLevelSets,
 thick] file {figuresAndData/dataForConditionalVarianceA0pt5B1pt5EqualAreas.tex};
\addplot[brown1, dashdotdotted, id= forLevelSets,
 thick] file {figuresAndData/dataForConditionalVarianceA0pt5B1pt5MeanZero.tex};
\addplot[ mark min,  
        blue5, id= forProb0pt2, thick] 
   file {figuresAndData/dataForConditionalVarianceA0pt5B1pt5Prob0pt2.tex};
\addplot[ mark min,  
        blue4, id= forProb0pt4, thick] 
   file {figuresAndData/dataForConditionalVarianceA0pt5B1pt5Prob0pt4.tex};
\addplot[ mark min,  
        blue3, id= forProb0pt6, thick] 
   file {figuresAndData/dataForConditionalVarianceA0pt5B1pt5Prob0pt6.tex};
\addplot[ mark min,  blue2,
         id= forProb0pt8, thick] 
   file {figuresAndData/dataForConditionalVarianceA0pt5B1pt5Prob0pt8.tex};
\end{axis}
\end{tikzpicture}
}
    \caption{The near-optimal performance of super-level intervals}
\label{fig:nearOptimalityOfSuperLevelSets}
\end{figure}

The variances incurred by the above families of silence intervals are depicted in Figure~\ref{fig:nearOptimalityOfSuperLevelSets}. Clearly super level interval are remarkably close to being optimal. Hence we can expect to get quite good approximations to the optimal silence interval by applying a couple of iterations of the centering algorithm, initializing it with a super-level interval.

\section{Bound on Rate distortion trade-offs for scalar unimodal densities, via Gauss inequality\label{section:rateDistortionBound}}
We now study the performance of probabilistic sampling for random variables that have symmetric, unimodal densities. Although log-concavity is not directly required for our result, that property is useful to preserve  the property of unimodality for the statistics of any random process with additive noise that is independent of the past and current states~\cite{ibragimov1956unimodal}. 

We use the Gauss inequality for the tails of scalar unimodal densities, to bound the rate-distortion curve of probabilistic sampling.
We consider silence intervals that are symmetric about a mode of the density. Recall that  the mode is a point such that the cumulative distribution function is convex everywhere to the left of the point, and is concave everywhere  to the right of this point. This point may be non-unique; nevertheless we pick a mode and denote it by~$\mu . $ 
We denote the silence set by~$ {\mathscr{S}}  = \left[ \mu - k , \mu + k \right] . $ We give upper bounds for: (i)~the sampling rate, and (ii)~the conditional variance given that the random variable falls within the silence interval.
\subsection{Upper bound on the sampling rate}
For a unimodal density~$p_X ( \cdot ) $, let~$\tau$ be defined as:
\begin{align*}
    \tau^2 & \triangleq 
    { \left( {\text{mean}} - {\text{mode}}   \right) }^2
    +    {\text{variance}} .
\end{align*}
The Gauss inequality~(
Section~1 in
~\cite{bickelKrieger1992extensionsOfChebyshev}
)
on the symmetric interval~$ {\mathscr{S}} $ around the mode, gives us:
\begin{gather} 
    {\mathbb{P}} \left[  \left\lvert X  - \mu  \right\rvert 
                           > k
                \right]  
     \  \le \
     \begin{cases}
         {\frac{4}{9}}
         {\frac{\tau^2}{k^2}} ,
         & {\text{if}} \ k \ge {\frac{2}{\sqrt{3}}} \tau ,
\\
         1 - {\frac{1}{\sqrt{3}}}
         {\frac{k}{\tau}} ,
         & {\text{if}}  \ 0 \le k \le {\frac{2}{\sqrt{3}}} \tau .
     \end{cases} 
     \label{eqn:gaussInequality}
\end{gather}
This gives us an upper bound on the rate  at which the state process~$X$ falls outside the silence interval, which is exactly the rate at which samples are generated.

\subsection{Upper bound on the mean-square error}
The mean-cum-mode~$ \mu $ is the least squares estimate of~$X$ given that it falls within the symmetric interval~$S.$
Because the density of~$X$ is unimodal, its conditional error variance of~$ X $ over the interval~${\mathscr{S}} $ 
can be bounded above by the variance of the uniform distribution over the same interval, as shown below. 

Let~$  {\mathscr{A}} $ be a measurable subset of~$  {\mathcal{S}} , 
$
and let $  {\mathbb{P}}_{U\left\lvert{\mathscr{S}}\right.} \left[ A \right] $ denote the probability mass of the  set~$ {\mathscr{A}} $ as per the uniform distribution over the interval~$ {\mathcal{S}} .
$
Consider the following two probabilities as functions of the positive real number~$t$:
\begin{align*}
    \Theta_{X\left\lvert{\mathscr{S}}\right.} \left( t \right) & \triangleq
  {\mathbb{P}}_{X\left\lvert{\mathscr{S}}\right.}
        \left[  \left\lvert X  - \mu  \right\rvert 
                           \le t  
                \right] ,
\\
    \Lambda_{U\left\lvert{\mathscr{S}}\right.} \left( t \right) & \triangleq 
    {\mathbb{P}}_{U\left\lvert{\mathscr{S}}\right.}
        \left[  \left\lvert X  - \mu  \right\rvert 
                         \le  t
                \right] .
\end{align*}
Clearly these two functions coincide at the extreme values of~$t$:
\begin{gather*}
  \Theta_{X\left\lvert{\mathscr{S}}\right.} \left( 0 \right) 
    \ =  \
    \Lambda_{U\left\lvert{\mathscr{S}}\right.} \left( 0 \right) \ = 0 , \ \ \text{and,} \\
  \Theta_{X\left\lvert{\mathscr{S}}\right.} \left( t \right) 
    \ =  \
    \Lambda_{U\left\lvert{\mathscr{S}}\right.} \left( t \right) \ = 1, \  \ \text{if}
    \  t \ge k .
\end{gather*}
On the interval~$
  \left[ 0 ,  k \right] 
 $ 
the function~$  
 \Lambda_{U\left\lvert{\mathscr{S}}\right.} \left( t \right) 
$ is linear and increasing.
On the same interval
the function~$  
 \Theta_{X\left\lvert{\mathscr{S}}\right.}  \left( t \right) 
$ is increasing. It is also concave. This is because
\begin{align*}
    \Theta_{X\left\lvert{\mathscr{S}}\right.} \left( t \right) & =
     {\mathbb{P}}_{X\left\lvert{\mathscr{S}}\right.}
        \left[  \mu \le   X \le    \mu + t
                \right] 
    +    {\mathbb{P}}_{X\left\lvert{\mathscr{S}}\right.}
        \left[   \mu  - t \le   X \le    \mu 
                \right] ,
\end{align*}
and both the terms on the right hand side are concave functions of~$t,$ since their derivatives with respective to~$t$ are non-increasing functions of~$t$~(because the density~$p_X$ is unimodal about the mode~$\mu$).

On the interval~$ \left[ 0 , k \right]$ 
the graph of~$  
  \Lambda_{U\left\lvert{\mathscr{S}}\right.}\left( t \right) 
$  is a straight line that intersects the graph of~$
\Theta_{X\left\lvert{\mathscr{S}}\right.} \left( t \right)
$ at the end points~$0 $ and~$ k .$ Hence the graph of~$
  \Lambda_{U\left\lvert{\mathscr{S}}\right.} \left( t \right) 
$ forms a chord for the graph of  the concave function~$
\Theta_{X\left\lvert{\mathscr{S}}\right.} \left( t \right)
$. Hence,  
conditioned on being in the interval~$
 {\mathscr{S}}
$ the density~$p_X ( \cdot ) $ is more peaked than the uniform density over the same interval, in the precise sense that:
\begin{align}
\Theta_{X\left\lvert{\mathscr{S}}\right.} \left( t \right)
    & \ge
    \Lambda_{U\left\lvert{\mathscr{S}}\right.} \left( t \right) , \ \text{for} \ 
    0 \ t \le k .
\end{align}
A direct consequence is that the corresponding variances obey:
\begin{align*}
    {\mathbb{E}}_{X\left\lvert{\mathscr{S}}\right.}
        \left[   {\left( X - \mu \right)}^2 
                 \left\lvert X \in {\mathscr{S}} \right.
                \right]  
    & \le
    {\mathbb{E}}_{U\left\lvert{\mathscr{S}}\right.} 
        \left[  
             {\left( X - \mu \right)}^2 
                 \left\lvert X \in {\mathscr{S}} \right.
                \right]  \ = \ {\frac{1}{3}} k^2 .
\end{align*}
This upper bound can be quite large if~$k$ is large. Therefore we cap it at  the variance~$\sigma^2$ to get the refined  upper bound below:
\begin{align}
    {\mathbb{E}}_{X\left\lvert{\mathscr{S}}\right.}
        \left[   {\left( X - \mu \right)}^2 
                 \left\lvert X \in {\mathscr{S}} \right.
                \right] 
    & \le
    \min{ \biggl\{  \sigma^2 , 
                    {\frac{1}{3}} k^2 
                    \biggr\}
                    } .
                \label{eqn:uniformDistributionVariance}
\end{align}
The mean squared error in estimating~$X,$ under probabilistic sampling that uses the silence set~$ {\mathcal{S}} $ is given by:
\begin{gather*}
  0 
    \times \left( 1 -  {\mathbb{P}}_X \left[  {\mathcal{S}} \right]  \right) 
  \ + \
    {\mathbb{E}}_{X\left\lvert{\mathscr{S}}\right.}
        \left[  {\left( X - \mu \right)}^2 
                 \left\lvert X \in {\mathscr{S}} \right.
                \right] 
    \times  {\mathbb{P}}_X \left[  {\mathcal{S}} \right]  .
\end{gather*}
Using Equation~\eqref{eqn:uniformDistributionVariance}, we can bound this quantity  above with:
\begin{gather}
     \begin{cases}
         {\frac{4}{9}}
         {\frac{\tau^2}{k^2}}
         \times
    \min{ \bigl\{  \sigma^2 , 
                    {\frac{1}{3}} k^2 
                    \bigr\}
                    } 
         ,
         & {\text{if}} \ k \ge {\frac{2}{\sqrt{3}}} \tau ,
\\
    \left(
         1 - {\frac{1}{\sqrt{3}}}
         {\frac{k}{\tau}} 
    \right)
         \times
    \min{ \bigl\{  \sigma^2 , 
                    {\frac{1}{3}} k^2 
                    \bigr\}
                    } 
         ,
         & {\text{if}}  \ 0 \le k \le {\frac{2}{\sqrt{3}}} \tau .
     \end{cases} 
     \label{eqn:boundOnDistortion}
\end{gather}
\subsection{A tight bound on the rate-distortion trade-off}
For any scalar unimodal density with a bounded second moment, and for a silence interval symmetric about the mode, we have shown that (1)~the sampling rate is bounded above by the RHS of~\eqref{eqn:gaussInequality}, and (2)~the squared error distortion is bounded above by the expression in~\eqref{eqn:boundOnDistortion}.
\begin{figure}
    \centering
\includegraphics[width=0.5\textwidth]{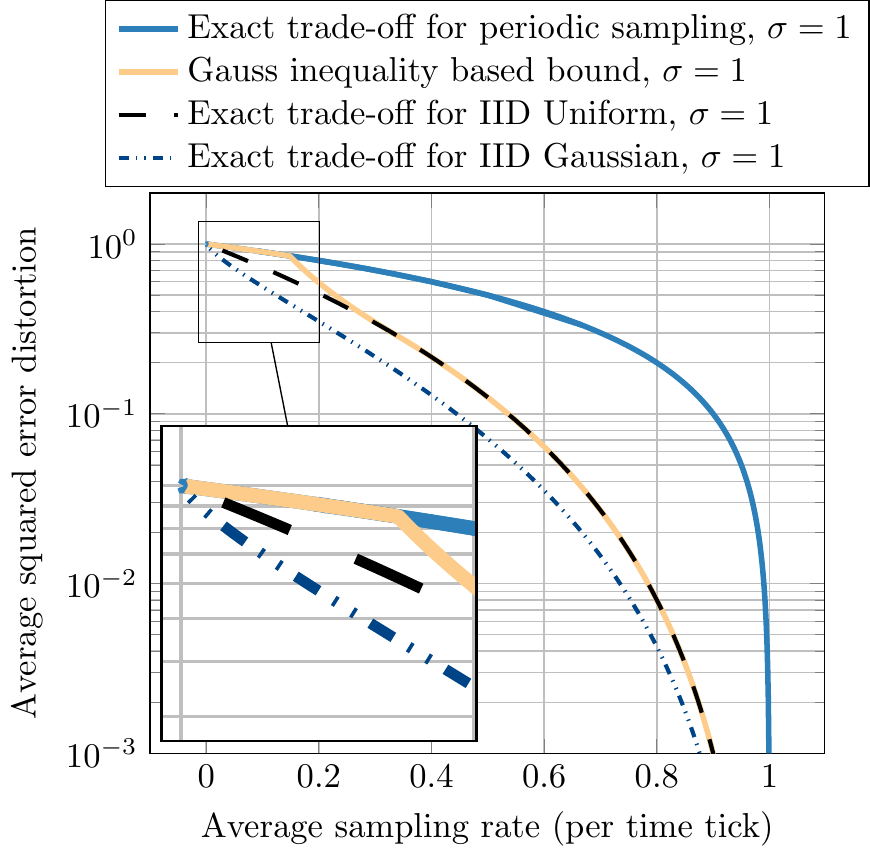}
\caption{Sampling rate versus Estimation error variance trade-off, for scalar IID proceses with  unimodal, symmetric densities.}
    \label{fig:gaussInequalityAndExactTradeOffs}
\end{figure}

We shall now apply these to the special case where the density is also symmetric about the mode. In this case, the quantity~$\tau^2$ found on the RHS of the Gauss inequality~\eqref{eqn:gaussInequality} simply equals the variance~$\sigma^2$. In Figure~\ref{fig:gaussInequalityAndExactTradeOffs} we plot our rate-distortion bound, and also the exact performances of periodic sampling, of optimal sampling when the density is uniform with~$\sigma = 1 $, and of optimal sampling when the density is Gaussian with~$\sigma = 1 .$

\section{Concluding remarks}
\subsection{Fast convergence of iterative centering}
\begin{figure}
    \centering
\includegraphics[width=0.5\textwidth]{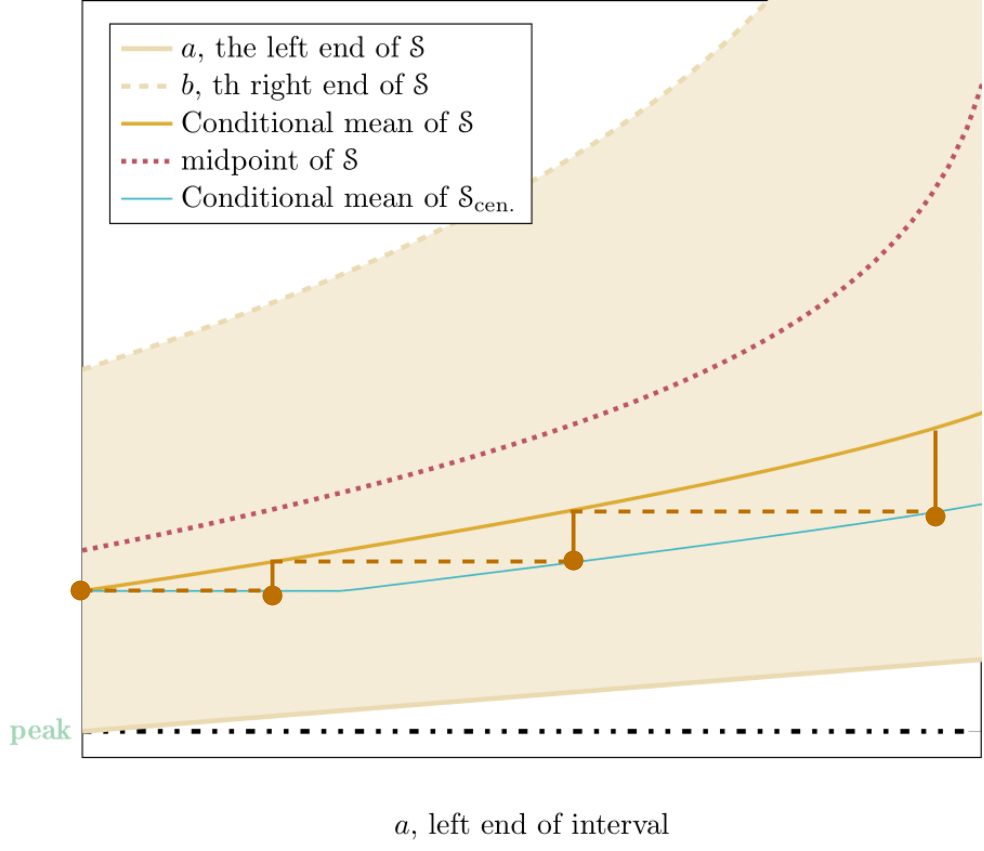}
\caption{Fast convergence of iterated centering}
    \label{fig:fastCOnvergenceOfCentering}
\end{figure}

For log-concave densities, the centering algorithm converges quite fast to optimum or near-optimum silence intervals, as illustrated in Figure~\ref{fig:fastCOnvergenceOfCentering}.
This figure concerns the exponential density, and   is a copy of Figure~\ref{subfig:conditionalMeansMeetAtLeftEnd}, but now showing the how quickly the conditional mean changes  under iterated  centering. The fast convergence happens basically because the conditional mean of the sliding interval  increases rather slowly as a function of the left end of the sliding interval. While each iteration  produces only a relatively modest change in the left end of the interval, it  nevertheless produces rather big changes in the value of the right end of the interval. This in turn produces big drops in the size of the interval, and also in its average distortion.

We conjecture that such fast convergence holds for every log-concave density. We think this because any tail of a log-concave density must be bounded above by an exponential decay. This shall ensure that the conditional mean of our  sliding interval changes relatively slowly, even as the midpoint and the length  change relatively fast. 

\subsection{Application to event-triggered sampling}
For even-triggered sampling of a random process,
a sub-optimal scheme
follows from our formulation for  individual random variables. But this needs to consider  forecasting and making  communication decisions over time. 
Firstly the sampled random process shall be correlated  over time, and secondly each sample is to be generated from within a horizon that  shall be  longer than just one tick. This  
entangles the design of the silence sets at any time tick to those at other time ticks.

For example consider a single sample to be generated from within a time horizon. We have to design a probability mass function~(PMF) of this sample being generated at any instant from within this horizon.
This PMF shall be non-uniorm, perhaps with small values at the start of the horizon, and  a peak at some suitable time instant such as the middle of horizon. At each instant,  we can propagate the nonlinear filter for  density of the process  at that instant, given that the sample has not been generated yet. We can then apply our centering algorithm on that density. 

\subsection{Utility of the Gauss bound}
The Gauss inequality is useful when the average sampling rate is in the range of approximately~0.3 to~0.9 samples per tick. It is in this range that we get the most benefits of interval based probabilistic sampling.
As Figure~\ref{fig:gaussInequalityAndExactTradeOffs} shows, there is almost no advantage to interval based sampling, if the average sampling rate is either close to zero or close to one tick per tick.
We can also see an order of magnitude drop in the average distortion,  if the sampling rate is between 0.6~and 0.9~samples per tick.

\section*{Acknowledgments}
M.~R. and V.~S. give thanks to Varun Bansal for help with programming in Python.
M.~R. also thanks the Automatic Control lab at KTH for hosting him during July 2016; this stay produced the main results recorded in this paper.
The authors thank Adam Molin for useful discussions on Example~\ref{example:circularDisk}.
The authors thank Chithrupa Ramesh, Takashi Tanaka, and Alexandre Proutiere for valuable discussions, and the anonymous reviewers of the 2023 American Control Conference, for identifying errors in an earlier version.

\bibliographystyle{IEEEtran}

\def\polhk#1{\setbox0=\hbox{#1}{\ooalign{\hidewidth
  \lower1.5ex\hbox{`}\hidewidth\crcr\unhbox0}}} \def\cprime{$'$}
\begin{thebibliography}{10}
\providecommand{\url}[1]{#1}
\csname url@samestyle\endcsname
\providecommand{\newblock}{\relax}
\providecommand{\bibinfo}[2]{#2}
\providecommand{\BIBentrySTDinterwordspacing}{\spaceskip=0pt\relax}
\providecommand{\BIBentryALTinterwordstretchfactor}{4}
\providecommand{\BIBentryALTinterwordspacing}{\spaceskip=\fontdimen2\font plus
\BIBentryALTinterwordstretchfactor\fontdimen3\font minus
  \fontdimen4\font\relax}
\providecommand{\BIBforeignlanguage}[2]{{%
\expandafter\ifx\csname l@#1\endcsname\relax
\typeout{** WARNING: IEEEtran.bst: No hyphenation pattern has been}%
\typeout{** loaded for the language `#1'. Using the pattern for}%
\typeout{** the default language instead.}%
\else
\language=\csname l@#1\endcsname
\fi
#2}}
\providecommand{\BIBdecl}{\relax}
\BIBdecl

\bibitem{miskowicz2015eventBasedControlBook}
M.~Miskowicz, Ed., \emph{Event-Based Control and Signal Processing}.\hskip 1em
  plus 0.5em minus 0.4em\relax CRC press, 2015.

\bibitem{lemmon2009sienaSurvey}
M.~Lemmon, ``Event-triggered feedback in control, estimation, and
  optimization,'' in \emph{Networked Control Systems}, 
        A.~Bemporad, M.~Heemels, and M.~Johansson,
  Eds.\hskip 1em plus 0.5em minus 0.4em\relax Springer London, 2010, vol. 406,
  pp. 293--358.

\bibitem{heemelsJohanssonTabuada2012}
W.~P. M.~H. Heemels, K.~H. Johansson, and P.~Tabuada, ``An introduction to
  event-triggered and self-triggered control,'' in \emph{2012 IEEE 51st IEEE
  Conference on Decision and Control (CDC)}, Dec 2012, pp. 3270--3285.

\bibitem{astrom-bernhardsson-2003}
K.~J. {\AA}str{\"o}m and B.~Bernhardsson, ``Systems with {L}ebesgue sampling,''
  in \emph{Directions in mathematical systems theory and optimization}, 
  0.4em\relax Berlin: Springer, 2003, vol. 286, pp. 1--13.

\bibitem{shiShiChen2015eventBasedEstimationBook}
D.~Shi, L.~Shi, and T.~Chen, \emph{Event-Based State Estimation}.\hskip 1em
  plus 0.5em minus 0.4em\relax Springer Cham, 2015.

\bibitem{rabiRameshJohansson2016siam}
\BIBentryALTinterwordspacing
M.~Rabi, C.~Ramesh, and K.~H. Johansson, ``Separated design of encoder and
  controller for networked linear quadratic optimal control,'' \emph{SIAM
  Journal on Control and Optimization}, vol.~54, no.~2, pp. 662--689, 2016.
\BIBentrySTDinterwordspacing

\bibitem{jhelumChakravorthyMahajan2020remoteEstimation}
J.~{Chakravorty} and A.~{Mahajan}, ``Remote estimation over a packet-drop
  channel with markovian state,'' \emph{IEEE Transactions on Automatic
  Control}, vol.~65, no.~5, pp. 2016--2031, 2020.

\bibitem{jhelumChakravorty2017fundamentalLimits}
J.~Chakravorty and A.~Mahajan, ``Fundamental limits of remote estimation of
  autoregressive markov processes under communication constraints,'' \emph{IEEE
  Transactions on Automatic Control}, vol.~62, no.~3, pp. 1109--1124, 2017.

\bibitem{khojasteh2020valueOfTimingInformation}
M.~J. Khojasteh, P.~Tallapragada, J.~Cortés, and M.~Franceschetti, ``The value
  of timing information in event-triggered control,'' \emph{IEEE Transactions
  on Automatic Control}, vol.~65, no.~3, pp. 925--940, 2020.

\bibitem{hajekPagingJournal}
\BIBentryALTinterwordspacing
B.~Hajek, K.~Mitzel, and S.~Yang, ``Paging and registration in cellular
  networks: jointly optimal policies and an iterative algorithm,'' \emph{IEEE
  Trans. Information Theory}, vol.~54, no.~2, pp. 608--622, 2008. [Online].
  Available: \url{http://dx.doi.org/10.1109/TIT.2007.913566}
\BIBentrySTDinterwordspacing

\bibitem{lipsaMartins2011}
G.~Lipsa and N.~Martins, ``Remote state estimation with communication costs for
  first-order lti systems,'' \emph{IEEE Trans. Automatic Control}, vol.~56,
  no.~9, pp. 2013--2025, 2011.

\bibitem{nayyarAndThreeProfessorsTransactions2013}
A.~Nayyar, T.~Ba{\c{s}}ar, D.~Teneketzis, and V.~Veeravalli, ``Optimal
  strategies for communication and remote estimation with an energy harvesting
  sensor,'' \emph{IEEE Trans. Automatic Control}, vol.~58, no.~9, pp.
  2246--2260, 2013.

\bibitem{andren2017cdcPaperWithNonConvexSilenceSets}
M.~T. {Andrén}, B.~{Bernhardsson}, A.~{Cervin}, and K.~{Soltesz}, ``On
  event-based sampling for lqg-optimal control,'' in \emph{2017 IEEE 56th
  Annual Conference on Decision and Control (CDC)}, Dec 2017, pp. 5438--5444.

\bibitem{molin2012unimodalBimodal}
\BIBentryALTinterwordspacing
A.~Molin and S.~Hirche, ``An iterative algorithm for optimal event-triggered
  estimation,'' \emph{IFAC Proceedings Volumes}, vol.~45, no.~9, pp. 64 -- 69,
  2012.
\BIBentrySTDinterwordspacing

\bibitem{molinHirche2017iterativeSilenceSets}
A.~{Molin} and S.~{Hirche}, ``Event-triggered state estimation: An iterative
  algorithm and optimality properties,'' \emph{IEEE Transactions on Automatic
  Control}, vol.~62, no.~11, pp. 5939--5946, 2017.

\bibitem{an1996logConcave}
\BIBentryALTinterwordspacing
M.~Y. An, ``{Log-concave Probability Distributions: Theory and Statistical
  Testing},'' Duke University, Tech. Rep., 1995. [Online]. Available:
  \url{https://econwpa.ub.uni-muenchen.de/econ-wp/game/papers/9611/9611002.pdf}
\BIBentrySTDinterwordspacing

\bibitem{henningssonAstrom2006mtnsLogConcave}
T.~Henningsson and K.~J. {\AA}str{\"o}m, ``Log-concave observers,'' in
  \emph{Proceedings of the 17th International Symposium on Mathematical Theory
  of Networks and Systems}, Kyoto, Japan, Jul. 2006.

\bibitem{trushkin1982}
A.~V. {T}rushkin, ``{S}ufficient conditions for {U}niqueness of a {L}ocally
  {O}ptimal {Q}uantizer for a class of {C}onvex error weighting functions,''
  \emph{{IEEE} Trans. {I}nformation theory}, vol. IT-28, no.~2, pp. 187--198,
  March 1982.

\bibitem{trushkin1984}
------, ``{M}onotony of {L}loyd's method {II} for log-concave density and
  convex error weighting function (corresp.),'' \emph{IEEE Trans. Information
  Theory}, vol.~30, no.~2, pp. 380--383, 1984.

\bibitem{kieffer1983uniquenessOfOptimalQuantizer}
J.~{Kieffer}, ``Uniqueness of locally optimal quantizer for log-concave density
  and convex error weighting function,'' \emph{IEEE Transactions on Information
  Theory}, vol.~29, no.~1, pp. 42--47, 1983.

\bibitem{ibragimov1956unimodal}
\BIBentryALTinterwordspacing
I.~Ibragimov, ``On the composition of unimodal distributions,'' \emph{Theor.
  Probability Appl.}, vol.~1, no.~2, pp. 255--260, 1956. 
\BIBentrySTDinterwordspacing

\bibitem{bickelKrieger1992extensionsOfChebyshev}
P.~J. Bickel and A.~M. Krieger, ``Extensions of {C}hebychev's inequality with
  applications,'' \emph{Probab. Math. Statist.}, vol.~13, no.~2, pp. 293--310
  (1993), 1992.

\end{thebibliography}
\def\polhk#1{\setbox0=\hbox{#1}{\ooalign{\hidewidth
  \lower1.5ex\hbox{`}\hidewidth\crcr\unhbox0}}} \def\cprime{$'$}

\end{document}